\documentclass[11pt,a4paper]{article}
\usepackage[margin=0.8in]{geometry}
\usepackage{amsmath,amssymb,amsthm,mathtools,mathrsfs}
\usepackage{xcolor}
\usepackage[colorlinks=true,linkcolor=blue,citecolor=blue,urlcolor=blue]{hyperref}
\usepackage{enumitem}
\usepackage{booktabs}
\allowdisplaybreaks
\numberwithin{equation}{section}

\newif\ifdraft \drafttrue

\theoremstyle{plain}
\newtheorem{theorem}{Theorem}[section]
\newtheorem{proposition}[theorem]{Proposition}
\newtheorem{lemma}[theorem]{Lemma}
\newtheorem{corollary}[theorem]{Corollary}
\theoremstyle{definition}
\newtheorem{definition}[theorem]{Definition}
\newtheorem{assumption}[theorem]{Assumption}
\theoremstyle{remark}
\newtheorem{remark}[theorem]{Remark}

\newcommand{\bO}{\mathcal{O}}
\newcommand{\bV}{\mathbb{V}}
\newcommand{\bH}{\mathbb{H}}
\newcommand{\bL}{\mathbb{L}}
\newcommand{\bD}{\mathbb{D}}
\newcommand{\bJ}{\mathbb{J}}
\newcommand{\bQ}{\mathbb{Q}}
\newcommand{\bE}{\mathbb{E}}
\newcommand{\bP}{\mathbb{P}}
\newcommand{\Nt}{N_{\tau}}
\newcommand{\vu}{\mathbf{u}}
\newcommand{\vd}{\mathbf{d}}
\newcommand{\vW}{\mathbf{W}}
\newcommand{\vze}{\boldsymbol{\zeta}}

\newcommand{\bR}{\mathbb{R}}

\newcommand{\Ah}{A_h}
\newcommand{\dmid}{\vd^{n+\frac12}}

\newcommand{\Ecal}{\mathcal{E}}
\newcommand{\Fcal}{\mathbb{F}}
\newcommand{\ip}[2]{\bigl(#1,#2\bigr)}
\newcommand{\iph}[2]{\bigl(#1,#2\bigr)_{h}}
\newcommand{\dual}[2]{\bigl\langle #1,#2\bigr\rangle}

\newcommand{\Isig}{\Theta}
\newcommand{\Gop}{\mathcal{G}}
\newcommand{\Ftil}{\widetilde F}
\newcommand{\ftil}{\widetilde f}
\newcommand{\Btil}{\widetilde B}
\newcommand{\btil}{\widetilde b}
\newcommand{\norm}[1]{\lVert #1\rVert}
\newcommand{\abs}[1]{\lvert #1\rvert}

\begin{document}


\title{\bf Energy-Compatible Finite Elements for Stochastic Nematic Liquid Crystals with Geometric Stratonovich Noise}

\author{
Akash Ashirbad Panda
\\[1mm]
{\footnotesize Department of Mathematics, Indian Institute of Technology Indore, Madhya Pradesh, India}
\\[-1mm]
{\small \texttt{akashpanda@iiti.ac.in}}
}


\maketitle

\begin{abstract}
We develop and analyse a structure-preserving finite element approximation of
the stochastic nematic liquid crystal system in
Ginzburg--Landau penalised form in two dimensions, with multiplicative It\^o noise in the
velocity and geometric Stratonovich noise in the director. The director noise
is treated by the implicit midpoint rule and the penalisation by a two-point
secant in $|\mathbf d|^{2}$.

A central result is a compatibility principle showing that, among continuous
field-independent two-point discretisations, exact potential-energy balance
and orthogonality to the midpoint noise uniquely determine the
$|\mathbf d|^{2}$-secant. Consequently, the stochastic part of the discrete
energy law contains no negative power of the penalty parameter. For the
isolated director-noise problem, the midpoint update is moreover an exact
isometry of strong order one.

For fixed penalty, we prove unconditional one-step solvability, localised
stability without an $h$--$\tau$ coupling, and convergence in probability to
a weak martingale solution. Under additional regularity, the coupled scheme
converges with rate
$O_{\mathbb P}(h+\tau^{1/2})$; for prescribed velocity the director is first
order in time, a behaviour also observed numerically in the coupled problem.
Spectral and finite element experiments confirm the structural identities and
the predicted convergence rates.
\end{abstract}

{\raggedright\noindent\textbf{Keywords.}
nematic liquid crystals; stochastic Ericksen--Leslie system;
finite element method; structure-preserving discretisation;
Ginzburg--Landau penalisation; Stratonovich noise;
weak martingale solutions.\par}


\section{Introduction}\label{sec:intro}

Nematic liquid crystals are fluids of rod-like molecules that retain a preferred
direction of alignment, described by a director field taking values on the unit sphere
\cite{deGennesProst1993}. The hydrodynamic theory of Ericksen
\cite{Ericksen1961} and Leslie \cite{Leslie1968} couples the Navier--Stokes
equations to an evolution equation for the director through an anisotropic elastic
stress. A fundamental feature of the model is its energy structure: the transport of
the director and the elastic stress balance exactly. To avoid the analytical
difficulties associated with the pointwise constraint $\abs{\vd}=1$, Lin and Liu
\cite{LinLiu1995} introduced a Ginzburg--Landau penalisation of strength
$\delta^{-2}$. The analysis and numerical approximation of the resulting penalised
system are by now well developed
\cite{LiuWalkington2000,BeckerFengProhl2008,Walkington2011}.

Thermal fluctuations, imperfections of the confining geometry, and fluctuations of
applied fields motivate stochastic perturbations of liquid-crystal models
\cite{DalbyMajumdarWuDond2024}. For the penalised stochastic nematic liquid-crystal
system, existence and uniqueness results are available
\cite{BHR2019,BMP2020}. The numerical problem is nevertheless delicate: a
discretisation must preserve not only the deterministic elastic--transport cancellation
and the chain rule associated with the Ginzburg--Landau energy, but also their
interaction with the geometric stochastic forcing. This interaction is the main subject
of the present work.

\subsection{The model}

Let $\bO\subset\bR^{2}$ be a bounded convex polygon and let $T>0$. We consider
\begin{align}
  d\vu + \bigl[(\vu\cdot\nabla)\vu-\Delta\vu+\nabla p\bigr]\,dt
  &=
  -\nabla\!\cdot\!\bigl(\nabla\vd\odot\nabla\vd\bigr)\,dt
  +\sigma(\vu)\,d\vW_{1},
  \label{eq:pde-u}\\
  \nabla\!\cdot\!\vu &=0,
  \label{eq:pde-div}\\
  d\vd+(\vu\cdot\nabla)\vd\,dt
  &=
  \left(\Delta\vd-\frac{1}{\delta^{2}}f(\vd)\right)dt
  +(\vd\times\vze)\circ d\vW_{2},
  \label{eq:pde-d}
\end{align}
in $(0,T]\times\bO$, with $\vu=0$ and $\partial_{n}\vd=0$ on $\partial\bO$, and initial data
$(\vu(0),\vd(0))=(\vu_{0},\vd_{0})$. Here
$\vu:[0,T]\times\bO\to\bR^{2}$ is the velocity,
$\vd:[0,T]\times\bO\to\bR^{3}$ the director, $p$ the pressure, and
$\delta>0$ the Ginzburg--Landau penalty parameter. We use the classical
potential
\begin{equation}\label{eq:f}
  f(\vd)=(\abs{\vd}^{2}-1)\vd,
  \qquad
  F(\vd)=\frac14\abs{\vd}^{4}-\frac12\abs{\vd}^{2},
  \qquad
  F'(\vd)[g]=f(\vd)\cdot g .
\end{equation}

The velocity equation is driven by multiplicative It\^o noise, whereas the
director equation contains the geometric Stratonovich noise
\[
G(\vd)=\vd\times\vze,
\qquad
G(\vd)\cdot\vd=0 .
\]
Thus, in the isolated noise equation, the stochastic forcing preserves
$\abs{\vd}$ exactly. This geometric property is consistent with the constrained
limit $\delta\to0$ \cite{BHR2019,BMP2019}. Throughout this paper, however,
$\delta>0$ is fixed: the director is $\bR^{3}$-valued and the constraint
$\abs{\vd}=1$ is not imposed.

\begin{remark}[Role of the geometric noise]\label{rem:constraint}
Although the penalised solution does not remain on the sphere, it is important that
the stochastic forcing itself does not create a spurious change in the director length.
The midpoint discretisation introduced below retains this property exactly for the
isolated noise problem; see Proposition \ref{prop:cayley} and
Remark \ref{rem:notpreserved}.
\end{remark}

\subsection{Related work and motivation}

Finite element methods for deterministic nematic liquid-crystal flows have a substantial
literature. Liu and Walkington \cite{LiuWalkington2000} proved convergence for the
penalised model, while Becker, Feng and Prohl \cite{BeckerFengProhl2008}
constructed energy-stable discretisations preserving the discrete analogue of the
Lin--Liu elastic--transport cancellation; see also
\cite{BarrettFengProhl2006,Walkington2011}. Constraint-preserving approximations of
sphere-valued sub-problems were studied in
\cite{BartelsProhl2007,BartelsLubichProhl2009}.

For stochastic nematic liquid crystals, weak martingale solutions and pathwise
uniqueness for the penalised system with multiplicative Gaussian noise were established
in \cite{BHR2019}; pure jump noise in Marcus canonical form was considered in
\cite{BMP2019}, and a large-deviation principle in \cite{BMP2020}. For the
sphere-constrained stochastic Ericksen--Leslie system, see
\cite{deBouardHocquetProhl2021,DuWang2022}. Related numerical analysis for SPDEs
includes finite element methods for stochastic Navier--Stokes equations
\cite{BrzezniakCarelliProhl2013,CarelliProhl2012,BreitDodgson2021,BreitProhl2023}
and structure-preserving approximation of the stochastic
Landau--Lifshitz--Gilbert equation
\cite{BanasBrzezniakNeklyudovProhl2014,BanasBrzezniakProhl2013}.

These deterministic numerical and stochastic analytical theories have largely developed
separately. Stochastic simulations of liquid-crystal models have been reported, for
example, in \cite{DalbyMajumdarWuDond2024}, but, to the best of our knowledge, a
convergence analysis of a finite element method for a stochastic nematic
liquid-crystal or Ericksen--Leslie system has not previously been established.

The main issue here is more structural than the construction of a convergent scheme
alone. The deterministic analysis relies simultaneously on the elastic--transport
cancellation and the chain rule for the penalisation. In the stochastic problem these
structures interact with the director noise. We therefore ask: \emph{How should geometric Stratonovich noise be discretised in the presence of a
penalised sphere constraint so that the discrete energy structure is preserved?}

Our answer is a compatibility principle: when the director noise is treated by the
implicit midpoint rule, exact reproduction of the potential-energy increment and
orthogonality to the midpoint noise determine the discrete nonlinearity uniquely,
within the natural class of continuous field-independent two-point maps. The resulting
nonlinearity is precisely the secant of the potential in $\abs{\vd}^{2}$.

\subsection{Main contributions}

The principal results of the paper are the following.
\begin{enumerate}[label=(C\arabic*),leftmargin=3em,itemsep=3pt]

\item \textbf{Structure-preserving finite element scheme and compatibility principle.}
We construct a fully discrete scheme using a pointwise divergence-free velocity space,
a mass-lumped $P_{1}$ director approximation, the implicit midpoint rule for the
geometric noise, and the $\abs{\vd}^{2}$-secant for the Ginzburg--Landau
nonlinearity. Theorem \ref{thm:characterisation} shows that, among continuous
field-independent two-point maps, this secant is uniquely characterised by exact
potential-energy balance and orthogonality to the midpoint noise for every applied
field. The principle extends beyond the quartic potential to general penalised sphere
constraints with geometric Stratonovich noise.

\item \textbf{Geometric exactness and the role of the midpoint rule.}
For the isolated director-noise equation, the midpoint update is the Cayley transform
of a skew-symmetric operator and is therefore an exact nodal isometry
(Proposition \ref{prop:cayley}). Since $\vW_{2}$ is scalar, the method has strong
order one without simulation of L\'evy areas. In contrast, an explicit treatment of the
noise produces an additional energy bias of order $O(\delta^{-2}\tau)$ and a
fluctuation of order $O(\delta^{-2}\tau^{1/2})$; uniform second-moment control as
$\delta\to0$ therefore requires $\tau=O(\delta^{4})$
(Proposition \ref{prop:negative}).

\item \textbf{Energy stability and convergence.}
The compatibility principle yields an exact discrete energy identity in which the
stochastic contribution contains no negative power of $\delta$
(Theorem \ref{thm:energy}). We prove unconditional solvability of each time step,
localised stability without any restriction coupling $\tau$ to $h$, and convergence
in probability to a weak martingale solution
(Theorems \ref{thm:stability} and \ref{thm:convergence}).

\item \textbf{Error estimates and temporal rates.}
Under additional spatial and temporal regularity assumptions, the fully coupled scheme
satisfies the rate $O_{\mathbb P}\!\left(h+\tau^{1/2}\right).$ The temporal half-order is caused by the multiplicative It\^o noise in the velocity
equation. For prescribed velocity, the director approximation attains first order in
time in both $\bL^{2}$ and $\bH^{1}$ under the corresponding regularity assumptions
(Proposition \ref{prop:drate}); the same behaviour is observed numerically for the
coupled problem.

\item \textbf{Numerical verification.}
Spectral experiments isolate the stochastic and algebraic structure, while genuine
finite element computations verify the spatial convergence rate, pointwise
divergence-free constraint, elastic--transport cancellation, and discrete energy law.
The experiments also confirm the predicted temporal rates and the advantage of the
midpoint treatment of the director noise.

\end{enumerate}

\medskip

\noindent
{\bf Organization of the paper.} \ Section \ref{sec:setting} introduces the functional setting and assumptions. Section \ref{sec:scheme} presents the finite element spaces and the fully discrete
scheme. Section \ref{sec:structure} establishes the compatibility principle and the
geometric properties of the director-noise discretisation. Section
\ref{sec:stability} derives the discrete energy law and stability estimates, while
Section \ref{sec:convergence} proves convergence to a weak martingale solution.
The error analysis is developed in Section \ref{sec:rates}, numerical experiments are
reported in Section \ref{sec:numerics}, and Section \ref{sec:remarks} concludes with
open problems.

\section{Functional setting}\label{sec:setting}

We collect the functional framework, structural identities, and assumptions used
throughout the analysis.

\subsection{Functional spaces and operators}

Let
\[
\mathcal{V}
:=
\left\{
\vu\in C_{c}^{\infty}(\bO;\bR^{2}):
\nabla\!\cdot\!\vu=0
\right\},
\]
and denote by $H$ and $\bV$ its closures in
$L^{2}(\bO;\bR^{2})$ and $H^{1}(\bO;\bR^{2})$, respectively. Then $\bV\hookrightarrow H\cong H'\hookrightarrow\bV'.$ For the director we write
\[
\bL^{2}=L^{2}(\bO;\bR^{3}),\qquad
\bH^{1}=H^{1}(\bO;\bR^{3}),\qquad
\bH^{2}=H^{2}(\bO;\bR^{3}).
\]
The notation $\ip{\cdot}{\cdot}$ and $\norm{\cdot}_{L^{2}}$ is used for the
$L^{2}(\bO)$ inner product and norm for scalar-, $\bR^{2}$- and
$\bR^{3}$-valued functions. Single bars $\abs{\cdot}$ denote Euclidean modulus
(or Lebesgue measure when applied to $\bO$). Since elements of $\bV$ vanish on
$\partial\bO$, the Poincar\'e inequality implies that $\norm{\nabla\cdot}_{L^{2}}$ is a
norm on $\bV$ equivalent to the $H^{1}$ norm.

For $\beta\in(0,1)$ and a Banach space $X$, we use, as in the compactness method of
Flandoli and Gatarek \cite{FlandoliGatarek1995}, the space
$W^{\beta,2}(0,T;X)$ with norm
\begin{equation}\label{eq:frac-def}
\norm{v}_{W^{\beta,2}(0,T;X)}^{2}
:=
\norm{v}_{L^{2}(0,T;X)}^{2}
+
\int_{0}^{T}\!\!\int_{0}^{T}
\frac{\norm{v(t)-v(s)}_{X}^{2}}
     {\abs{t-s}^{1+2\beta}}
\,dt\,ds .
\end{equation}

Let $\mathscr{A}$ denote the Stokes operator,
\[
\dual{\mathscr{A}\vu}{v}
=
\ip{\nabla\vu}{\nabla v},
\qquad
D(\mathscr{A})
=
\bV\cap H^{2}(\bO;\bR^{2}),
\]
and let
\[
A=-\Delta,
\qquad
D(A)=\bigl\{\vd\in\bH^{2}:\partial_{n}\vd=0\text{ on }\partial\bO\bigr\}.
\]
Convexity of $\bO$ supplies the $H^{2}$ regularity of the Stokes and Neumann problems
behind both identifications \cite{KelloggOsborn1976,Grisvard1985}. In two dimensions,
\begin{equation}\label{eq:domemb}
D(\mathscr{A})
\hookrightarrow
L^{\infty}\cap W^{1,4},
\qquad
D(A)\hookrightarrow L^{\infty}.
\end{equation}
These embeddings will be used in the time-regularity estimates of
\S\ref{sec:convergence}; see Remark \ref{rem:dual}.

For the nonlinear terms define
\begin{align}
b(\vu,v,z)
&:=
\sum_{i,j}
\int_{\bO}
\vu^{(i)}\partial_{x_i}v^{(j)}z^{(j)}\,dx,
&
\dual{B(\vu,v)}{z}
&:=b(\vu,v,z),
\label{eq:b}
\\
\btil(\vu,\vd,e)
&:=
\sum_i
\int_{\bO}
\vu^{(i)}\partial_{x_i}\vd\cdot e\,dx,
&
\dual{\Btil(\vu,\vd)}{e}
&:=\btil(\vu,\vd,e),
\label{eq:btil}
\\
m(\vd_1,\vd_2,\vu)
&:=
-\sum_{i,j,k}
\int_{\bO}
\partial_{x_j}\vd_1^{(k)}
\partial_{x_i}\vd_2^{(k)}
\partial_{x_j}\vu^{(i)}\,dx,
&
\dual{M(\vd_1,\vd_2)}{\vu}
&:=m(\vd_1,\vd_2,\vu).
\label{eq:m}
\end{align}
We write $B(\vu):=B(\vu,\vu)$ and $M(\vd):=M(\vd,\vd)$. Thus
\[
\dual{M(\vd)}{\vu}
=
-\int_{\bO}
(\nabla\vd\odot\nabla\vd):\nabla\vu\,dx,
\]
which is the weak form of the elastic stress in \eqref{eq:pde-u}.

\subsection{Structural identities and assumptions}

The deterministic energy structure is encoded in the following identities.

\begin{lemma}\label{lem:cancel}
For $\vu\in\bV$ and $\vd\in D(A)$,
\begin{equation}\label{eq:cont-cancel}
b(\vu,v,v)=0,\qquad
\btil(\vu,\vd,\vd)=0,\qquad
m(\vd,\vd,\vu)
=
-\btil(\vu,\vd,A\vd),
\qquad
\btil(\vu,\vd,f(\vd))=0.
\end{equation}
\end{lemma}

\begin{proof}
The first two identities follow from $\nabla\!\cdot\!\vu=0$ and $\vu|_{\partial\bO}=0$.
Moreover,
\[
-\btil(\vu,\vd,A\vd)
=
\int_{\bO}
\Delta\vd\cdot(\vu\cdot\nabla)\vd\,dx
=
-\int_{\bO}
(\nabla\vd\odot\nabla\vd):\nabla\vu\,dx,
\]
where the boundary terms vanish because $\vu=0$ and $\partial_{n}\vd=0$ on
$\partial\bO$, and so does
$\frac12\int_{\bO}(\nabla\!\cdot\!\vu)\abs{\nabla\vd}^{2}\,dx$. Finally, since
$f(\vd)=\ftil(\abs{\vd}^{2})\vd$,
\[
\btil(\vu,\vd,f(\vd))
=
\int_{\bO}
\vu\cdot\nabla F(\vd)\,dx
=0 .
\]
\end{proof}

\begin{remark}[Sign convention]\label{rem:sign}
The minus sign in the third identity of \eqref{eq:cont-cancel} is essential:
it gives the cancellation between the elastic stress and director transport
in the energy balance and determines the sign of the discrete elastic operator
\eqref{eq:Mh}.
\end{remark}

We impose the following assumptions.

\begin{assumption}[Noise]\label{ass:noise}
$\vW_{1}$ is an $H$-valued Wiener process with reproducing kernel Hilbert space
$H_{0}\hookrightarrow H$, the embedding being Hilbert--Schmidt, and $\vW_{2}$
is a real-valued Wiener process independent of $\vW_{1}$. Both are defined on
a filtered probability space
$(\Omega,\mathcal F,\Fcal,\bP)$ satisfying the usual conditions.

The mapping
\[
\sigma:H\longrightarrow L_{2}(H_{0},H)
\]
satisfies
\[
\norm{\sigma(\vu)}_{\mathrm{HS}}^{2}
\le
K_{\sigma}\bigl(1+\norm{\vu}_{L^{2}}^{2}\bigr),
\qquad
\norm{\sigma(\vu_{1})-\sigma(\vu_{2})}_{\mathrm{HS}}^{2}
\le
L_{\sigma}\norm{\vu_{1}-\vu_{2}}_{L^{2}}^{2}.
\]
\end{assumption}

\begin{assumption}[Applied field]\label{ass:h}
The applied field satisfies
\[
\vze\in W^{1,\infty}(\bO;\bR^{3}).
\]
We set
\[
G(\vd)=\vd\times\vze,
\qquad
\Gop(x)=-[\vze(x)]_{\times},
\]
so that
$G(\vd)(x)=\Gop(x)\vd(x)$ and $G^{2}=G\circ G$.
\end{assumption}

\begin{assumption}[Potential]\label{ass:F}
We use the classical Ginzburg--Landau potential
\[
\Ftil(\varsigma)
=
\frac12\varsigma^{2}-\varsigma,
\qquad
\ftil(\varsigma)
=
\Ftil'(\varsigma)
=
\varsigma-1,
\]
and therefore
\[
f(\vd)
=
(\abs{\vd}^{2}-1)\vd,
\qquad
F(\vd)
=
\frac12\Ftil(\abs{\vd}^{2}),
\qquad
F'(\vd)[g]
=
f(\vd)\cdot g.
\]
\end{assumption}

\begin{remark}[Why the quartic, and where generality survives]
\label{rem:quartic}
The analysis in \S\S\ref{sec:stability}--\ref{sec:rates} is restricted to the
classical quartic because the uniform one-sided Lipschitz estimate
\eqref{eq:osl}, required in Proposition \ref{prop:solvable}, is established
here only in this case. The structural compatibility principle is more
general: Theorem \ref{thm:characterisation} requires only
$\Ftil\in C^{1}$, while Lemma \ref{lem:fstar} applies to polynomial
$\Ftil$.
\end{remark}

\begin{remark}[Normalisation]\label{rem:normal}
Adding a constant to $\Ftil$ changes neither $f$ nor the divided difference
$\kappa$ of Definition \ref{def:fstar}; it only shifts the energy by a
constant. Thus the normalisation
$\Ftil(\varsigma)=\frac12\varsigma^{2}-\varsigma$ entails no loss of
generality for the scheme.
\end{remark}

\subsection{Energy and coercivity}

Define the director and total energies by
\begin{equation}\label{eq:energy}
\Phi(\vd)
:=
\frac12\norm{\nabla\vd}_{L^{2}}^{2}
+
\frac{1}{\delta^{2}}
\int_{\bO}F(\vd)\,dx,
\qquad
\Ecal(\vu,\vd)
:=
\frac12\norm{\vu}_{L^{2}}^{2}
+
\Phi(\vd).
\end{equation}
For $(\vu_{0},\vd_{0})\in H\times\bH^{1}$ these quantities are finite, since
in two dimensions $\bH^{1}\hookrightarrow L^{q}$ for every $q<\infty$.

The quartic potential also provides the coercivity needed below.

\begin{lemma}[Coercivity of the potential]\label{lem:coerc}
There exist $c_{F},C_{F}>0$ such that
\[
F(\vd)
\ge
c_{F}\abs{\vd}^{4}-C_{F},
\qquad
\vd\in\bR^{3}.
\]
Consequently,
\begin{equation}\label{eq:coerc}
\Phi(\vd)
\ge
\frac12\norm{\nabla\vd}_{L^{2}}^{2}
+
\frac{c_{F}}{\delta^{2}}
\norm{\vd}_{L^{4}}^{4}
-
\frac{C_{F}\abs{\bO}}{\delta^{2}}.
\end{equation}
In particular, $\Ecal$ is bounded below by
$-C_{F}\abs{\bO}\delta^{-2}$ and controls
$\norm{\vd}_{L^{4}}$.
\end{lemma}

\begin{proof}
Since
\[
\Ftil(\varsigma)
=
\frac12\varsigma^{2}-\varsigma
\ge
\frac14\varsigma^{2}-1,
\qquad
\varsigma\ge0,
\]
and $F=\frac12\Ftil(\abs{\vd}^{2})$, one may take
$c_{F}=1/8$ and $C_{F}=1/2$. Integration over $\bO$ gives
\eqref{eq:coerc}.
\end{proof}

In the Gr\"onwall estimates below we therefore work, without further comment,
with the nonnegative shifted energy $\Ecal_{h}+C_{F}\abs{\bO}\delta^{-2}.$ The two-dimensional interpolation inequalities, Aubin--Lions compactness,
discrete Gr\"onwall inequality, and Burkholder--Davis--Gundy inequality used
throughout are collected in Appendix \ref{app:tools}.

\subsection{Abstract formulation and well-posedness}

The derivative of the director energy satisfies
\begin{equation}\label{eq:Phi}
\Phi'(\vd)[g]
=
\dual{
A\vd+\frac{1}{\delta^{2}}f(\vd)
}{g}.
\end{equation}
Hence \eqref{eq:pde-u}--\eqref{eq:pde-d} can be written as
\begin{align}
d\vu
+
\bigl[
\mathscr{A}\vu+B(\vu)+M(\vd)
\bigr]dt
&=
\sigma(\vu)\,d\vW_{1},
\label{eq:abs-u}
\\
d\vd
+
\left[
A\vd+\Btil(\vu,\vd)
+\frac{1}{\delta^{2}}f(\vd)
\right]dt
&=
G(\vd)\circ d\vW_{2}.
\label{eq:abs-d}
\end{align}

\begin{definition}[Weak martingale solution]\label{def:wms}
A weak martingale solution of
\eqref{eq:abs-u}--\eqref{eq:abs-d} is a system
\[
(\Omega,\mathcal F,\Fcal,\bP,\vW,\vu,\vd),
\qquad
\vW=(\vW_{1},\vW_{2}),
\]
such that $(\vu,\vd)$ is progressively measurable, has $\bP$-a.s.\
continuous paths in $H\times\bH^{1}$,
\[
\bE\sup_{t\in[0,T]}
\left(
\norm{\vu(t)}_{L^{2}}^{2}
+
\norm{\nabla\vd(t)}_{L^{2}}^{2}
\right)
+
\bE\int_{0}^{T}
\left(
\norm{\nabla\vu}_{L^{2}}^{2}
+
\norm{A\vd}_{L^{2}}^{2}
\right)dt
<\infty,
\]
and satisfies the variational form of
\eqref{eq:abs-u}--\eqref{eq:abs-d} $\bP$-a.s.\ against test functions in
$\bV\times\bL^{2}$, with
\[
G(\vd)\circ d\vW_{2}
=
\frac12G^{2}(\vd)\,dt
+
G(\vd)\,d\vW_{2}.
\]
\end{definition}

\begin{theorem}[{\cite[Thm.~3.1]{BHR2019}}]\label{thm:existence}
Under Assumptions \ref{ass:noise}--\ref{ass:F}, for
$(\vu_{0},\vd_{0})\in H\times\bH^{1}$,
system \eqref{eq:abs-u}--\eqref{eq:abs-d} possesses a weak martingale
solution in the sense of Definition \ref{def:wms}, and the solution is
pathwise unique in two dimensions.
\end{theorem}

\begin{remark}[Regularity of the applied field]\label{rem:field-reg}
Assumption \ref{ass:h} strengthens the
$W^{1,3}\cap L^{\infty}$ condition used in
\cite{BHR2019,BMP2020} to $W^{1,\infty}$. This is needed only in the
estimate of $\mathcal N^{n}$ in \eqref{eq:Nn} and in the product
interpolation estimate of Lemma \ref{lem:disc}\emph{(vi)}. The first can also
be obtained under $W^{1,3}\cap L^{\infty}$ using
\[
\norm{\nabla I_{h}(v\times\vze)}_{L^{2}}
\le
C\norm{\nabla\vze}_{L^{3}}\norm{v}_{L^{6}}
+
C\norm{\vze}_{L^{\infty}}\norm{\nabla v}_{L^{2}},
\]
but retaining $W^{1,\infty}$ avoids an additional elementwise argument and is
sufficient for the applications considered here.
\end{remark}

\section{The finite element scheme}\label{sec:scheme}

We now construct the fully discrete approximation. Three ingredients are essential:
a pointwise divergence-free velocity space, a mass-lumped $P_{1}$ approximation of
the director, and the joint treatment of the geometric noise and the
Ginzburg--Landau potential. The latter two will be shown in
\S\ref{sec:structure} to satisfy a rigid compatibility principle. Throughout this section, Assumptions \ref{ass:noise}--\ref{ass:F} are in force.

\subsection{Discrete spaces and basic properties}\label{sub:spaces}

\begin{assumption}[Discretisation]\label{ass:disc}
Let $\{\mathcal T_h\}_{h\in(0,h_0]}$ be a shape-regular and quasi-uniform
family of triangulations of $\bO$, with vertices
$\{x_j\}_{j=1}^{J_h}$.

\begin{enumerate}[label=\emph{(D\arabic*)},leftmargin=3.2em,itemsep=2pt]

\item
$(\bV_h,\bQ_h)\subset H_0^1(\bO;\bR^2)\times L_0^2(\bO)$, with
$\bQ_h:=\nabla\!\cdot\!\bV_h$, is an inf--sup
stable velocity--pressure pair \cite{GiraultRaviart1986}:
\begin{equation}\label{eq:infsup}
 \inf_{q_h\in\bQ_h\setminus\{0\}}
 \sup_{v_h\in\bV_h\setminus\{0\}}
 \frac{\ip{\nabla\!\cdot v_h}{q_h}}
      {\norm{\nabla v_h}_{L^2}\norm{q_h}_{L^2}}
 \ge \beta_\star>0 ,
\end{equation}
uniformly in $h$. Setting
\[
\bJ_h
:=
\left\{
v_h\in\bV_h:
\ip{\nabla\!\cdot v_h}{q_h}=0
\quad\forall q_h\in\bQ_h
\right\},
\]
we further require
\begin{equation}\label{eq:pwdf}
\nabla\!\cdot v_h=0
\quad\text{pointwise in }\bO,
\qquad v_h\in\bJ_h,
\end{equation}
and $\inf_{w_h\in\bJ_h}
\norm{\nabla(v-w_h)}_{L^2}
\longrightarrow 0$ \ for every $v\in\bV$\,.
Scott--Vogelius spaces of sufficiently high degree on meshes without
singular vertices provide a standard example
\cite{ScottVogelius1985}.

\item
The director space is the continuous piecewise affine space
\[
\bD_h
:=
\left\{
e_h\in C(\bO;\bR^3):
e_h|_{\mathcal K}\in P_1(\mathcal K)^3
\quad\forall\mathcal K\in\mathcal T_h
\right\},
\]
with nodal basis $\{\varphi_j\}$ and nodal interpolant $I_h$.

\item
The triangulation is weakly acute: \ $\ip{\nabla\varphi_i}{\nabla\varphi_j}\le0, \ i\ne j.$ Writing $\omega_{ij} :=
-\ip{\nabla\varphi_i}{\nabla\varphi_j}\ge0,$ we have
\begin{equation}\label{eq:edge}
\iph{\Ah\phi_h}{\psi_h}
=
\sum_{i<j}\omega_{ij}
\bigl(\phi_h(x_i)-\phi_h(x_j)\bigr)
\cdot
\bigl(\psi_h(x_i)-\psi_h(x_j)\bigr).
\end{equation}

\item
The lumped inner product is
\[
\iph{\phi}{\psi}
:=
\sum_{j=1}^{J_h}
\varpi_j\,\phi(x_j)\cdot\psi(x_j),
\qquad
\varpi_j:=\int_{\bO}\varphi_j\,dx>0,
\]
with associated norm $\norm{\cdot}_h$.

\end{enumerate}
\end{assumption}

The discrete Laplacian $\Ah:\bD_h\to\bD_h$ is defined by
\begin{equation}\label{eq:Ah}
\iph{\Ah\phi_h}{e_h}
=
\ip{\nabla\phi_h}{\nabla e_h},
\qquad
\phi_h,e_h\in\bD_h .
\end{equation}

We use the $L^2$-projections
\[
P_h:H\to\bJ_h,
\qquad
\Pi_h:L_0^2(\bO)\to\bQ_h,
\]
the Ritz projection $R_h:\bH^1\to\bD_h$ defined by
\begin{equation}\label{eq:ritz}
\ip{\nabla R_h\vd}{\nabla e_h}
+\ip{R_h\vd}{e_h}
=
\ip{\nabla\vd}{\nabla e_h}
+\ip{\vd}{e_h},
\qquad e_h\in\bD_h,
\end{equation}
and the Stokes projection $S_h:\bV\to\bJ_h$,
\begin{equation}\label{eq:stokesproj}
\ip{\nabla S_hv}{\nabla w_h}
=
\ip{\nabla v}{\nabla w_h},
\qquad w_h\in\bJ_h .
\end{equation}
Consequently,
\begin{equation}\label{eq:stokesproj-prop}
\norm{\nabla S_hv}_{L^2}
\le\norm{\nabla v}_{L^2},
\qquad
\norm{\nabla(v-S_hv)}_{L^2}
\le
\inf_{w_h\in\bJ_h}
\norm{\nabla(v-w_h)}_{L^2}
\longrightarrow0 .
\end{equation}

\begin{remark}[The domain and the mesh]\label{rem:domain}
Convexity of $\bO$ gives the $H^{2}$ regularity behind $D(A)$ and $D(\mathscr{A})$, and
polygonality lets affine triangulations represent $\bO$ exactly. Assumption
\ref{ass:disc} is nonempty. On a rectangle, the uniform triangulation with all diagonals in
one direction is weakly acute --- its off-diagonal stiffness entries are $-1$ and
$-\tfrac12$ --- and has no interior singular vertices, so weak acuteness and a pointwise
divergence-free Scott--Vogelius pair can be imposed together. The pressure space must then
be $\bQ_h=\nabla\!\cdot\!\bV_h$ rather than all of $P_{k-1}^{\rm disc}\cap L^2_0$: with
Dirichlet velocity the divergence misses part of the latter near the corners, and for $k=4$
on a square we find a rank deficiency of two beyond the constants, independent of $h$. The
velocity degree exceeds what the spatial rate uses, that rate being limited by the $P_1$
director; see \S\ref{sec:remarks}.
\end{remark}

The following properties will be used repeatedly.

\begin{lemma}[Discrete estimates]\label{lem:disc}
Under Assumption \ref{ass:disc}, with constants independent of $h$:

\begin{enumerate}[label=\emph{(\roman*)},leftmargin=2.4em,itemsep=2pt]

\item
\[
c_1\norm{e_h}_{L^2}
\le\norm{e_h}_h
\le c_2\norm{e_h}_{L^2}.
\]

\item
$\Ah$ is symmetric and positive semidefinite with respect to
$\iph{\cdot}{\cdot}$.

\item
\[
\norm{\nabla e_h}_{L^2}
\le Ch^{-1}\norm{e_h}_{L^2},
\qquad
\norm{\Ah e_h}_h
\le Ch^{-1}\norm{\nabla e_h}_{L^2},
\qquad
\norm{e_h}_{L^\infty}
\le Ch^{-1}\norm{e_h}_{L^2},
\]
and the analogous inverse estimates hold on $\bV_h$.

\item
The lumped product is nodal. Hence pointwise algebraic identities in
$\bR^3$ are inherited nodewise by the corresponding lumped pairings.

\item
For every $q\in[1,\infty)$,
\begin{equation}\label{eq:lumpLq}
c_q\norm{e_h}_{L^q}^q
\le
\sum_j\varpi_j\abs{e_h(x_j)}^q
\le
C_q\norm{e_h}_{L^q}^q .
\end{equation}

\item
For $\vze\in W^{1,\infty}(\bO;\bR^3)$,
\begin{equation}\label{eq:interp}
\norm{\nabla I_h(e_h\times\vze)}_{L^2}
\le
C\norm{\vze}_{W^{1,\infty}}
\norm{e_h}_{H^1}.
\end{equation}

\item
For $\phi_h\in\bD_h$,
\begin{equation}\label{eq:dGN}
\norm{\nabla\phi_h}_{L^4}
\le
C\norm{\nabla\phi_h}_{L^2}^{1/2}
\bigl(
\norm{\nabla\phi_h}_{L^2}
+\norm{\Ah\phi_h}_h
\bigr)^{1/2}.
\end{equation}

\end{enumerate}
\end{lemma}

\begin{proof}
Items (i)--(iv) are standard. Item (v) follows elementwise from equivalence of
norms on $P_1(\mathcal K)^3$ and scaling to the reference triangle.

For (vi), put $g=e_h\times\vze$. On each $\mathcal K$,
\[
\abs{\nabla I_hg}_{\mathcal K}
\le
C\left(
\norm{\vze}_{L^\infty}\abs{\nabla e_h}_{\mathcal K}
+
\norm{\nabla\vze}_{L^\infty}
\norm{e_h}_{L^\infty(\mathcal K)}
\right).
\]
Squaring, summing over $\mathcal T_h$, and using local norm equivalence gives
\eqref{eq:interp}. Estimate \eqref{eq:dGN} is the standard discrete
Gagliardo--Nirenberg estimate; see, e.g.\
\cite[Ch.~1]{BanasBrzezniakNeklyudovProhl2014}.
\end{proof}


\begin{remark}[Why $P_1$ with mass lumping]\label{rem:whyP1}
The compatibility identities of \S\ref{sec:structure} are pointwise and do
not by themselves require mass lumping. Mass lumping is needed for the
stronger geometric property of Proposition \ref{prop:cayley}: the isolated
noise step decouples node by node and therefore preserves
$\abs{\vd(x_j)}$ exactly. This motivates the mass-lumped $P_1$ director
space used throughout.
\end{remark}

\subsection{Discrete operators and the potential term}

For the velocity convection we use the skew-symmetric form
\[
\widehat b(w;v,z)
:=
\frac12\bigl[b(w,v,z)-b(w,z,v)\bigr],
\qquad
\widehat b(w;v,v)=0.
\]

The discrete elastic operator is defined by
\begin{equation}\label{eq:Mh}
\dual{M_h(\vd_h)}{v_h}
:=
-\btil(v_h,\vd_h,\Ah\vd_h),
\qquad
\vd_h\in\bD_h,\quad v_h\in\bV_h.
\end{equation}
Thus the continuous elastic--transport cancellation is reproduced exactly:
\begin{equation}\label{eq:BFP}
\dual{M_h(\vd_h)}{v_h}
+
\btil(v_h,\vd_h,\Ah\vd_h)
=0.
\end{equation}

For $\xi,\eta\in\bR^3$, set $\varrho:=\frac12(\xi+\eta)$ and define the $|\vd|^2$-secant by
\begin{equation}\label{eq:fstar}
f_\star(\xi,\eta)
:=
\kappa(\xi,\eta)\varrho,
\qquad
\kappa(\xi,\eta)
:=
\frac{
\Ftil(\abs{\eta}^2)-\Ftil(\abs{\xi}^2)
}{
\abs{\eta}^2-\abs{\xi}^2
},
\end{equation}
with the quotient understood by its continuous extension when
$\abs{\xi}=\abs{\eta}$.

\begin{definition}[$|\vd|^2$-secant]\label{def:fstar}
For the quartic Ginzburg--Landau potential,
\begin{equation}\label{eq:kappa-closed}
\kappa(\xi,\eta)
=
\frac12\bigl(\abs{\xi}^2+\abs{\eta}^2\bigr)-1,
\end{equation}
and hence
\[
f_\star(\xi,\eta)
=
\frac14
\bigl(\abs{\xi}^2+\abs{\eta}^2-2\bigr)
(\xi+\eta).
\]
For finite element functions, $f_\star$ is understood nodewise, equivalently
through its nodal interpolant in $\bD_h$.
\end{definition}

The elementary estimates required below are collected next.

\begin{lemma}[Properties of $f_\star$]\label{lem:fstarprop}
For $\xi,\eta,\eta_1,\eta_2\in\bR^3$:

\begin{enumerate}[label=\emph{(\roman*)},leftmargin=2.4em,itemsep=2pt]

\item
$f_\star$ is a polynomial map and $f_\star(\xi,\xi)=f(\xi).$

\item
\[
\abs{f_\star(\xi,\eta)-f(\eta)}
\le
C(1+\abs{\xi}+\abs{\eta})^2\abs{\eta-\xi}.
\]

\item
\[
\abs{f_\star(\xi,\eta)}
\le
C(1+\abs{\xi}+\abs{\eta})^3.
\]

\item
\[
f_\star(\xi,\eta)\cdot\eta
\ge
\frac18\abs{\eta}^4
-
C(1+\abs{\xi})^4.
\]

\item
With $z=\eta_1-\eta_2$,
\begin{equation}\label{eq:osl}
\bigl(
f_\star(\xi,\eta_1)-f_\star(\xi,\eta_2)
\bigr)\cdot z
\ge
-\frac12\abs{z}^2.
\end{equation}

\end{enumerate}
\end{lemma}

\begin{proof}
Items (i)--(iii) follow immediately from \eqref{eq:kappa-closed}.
For (iv),
\[
f_\star(\xi,\eta)\cdot\eta
=
\frac14\abs{\eta}^4
+\frac14\abs{\xi}^2\abs{\eta}^2
-\frac12\abs{\eta}^2
+
\frac14
\bigl(\abs{\xi}^2+\abs{\eta}^2-2\bigr)\xi\cdot\eta,
\]
and Young's inequality gives the assertion.

For (v), set
\[
\mu:=\frac12(\eta_1+\eta_2),
\qquad
\kappa_i:=\kappa(\xi,\eta_i).
\]
Since $\kappa_1-\kappa_2=\mu\cdot z$,
\begin{equation}\label{eq:oslsplit}
\begin{aligned}
&
\bigl(
f_\star(\xi,\eta_1)-f_\star(\xi,\eta_2)
\bigr)\cdot z =
\frac{\kappa_1+\kappa_2}{4}\abs{z}^2
+
\frac12(\mu\cdot z)(\xi\cdot z)
+
\frac12(\mu\cdot z)^2 .
\end{aligned}
\end{equation}
The last two terms are bounded below by
$-\frac18(\xi\cdot z)^2$, while
$\kappa_1+\kappa_2\ge\abs{\xi}^2-2$. Therefore
\[
\bigl(
f_\star(\xi,\eta_1)-f_\star(\xi,\eta_2)
\bigr)\cdot z
\ge
\left(
\frac18\abs{\xi}^2-\frac12
\right)\abs{z}^2
\ge
-\frac12\abs{z}^2.
\]
\end{proof}

\begin{remark}[Beyond the quartic]\label{rem:osl-general}
The secant construction \eqref{eq:fstar} and the structural identities of
\S\ref{sec:structure} extend to more general potentials. The uniform
one-sided estimate \eqref{eq:osl}, however, is established here only for the
classical quartic, which is why the stability and error analysis is formulated
under Assumption \ref{ass:F}.
\end{remark}

\subsection{The fully discrete scheme}

Let
\[
\tau=\frac{T}{\Nt},
\qquad
t_n=n\tau,
\qquad
a_n:=\vW_2(t_{n+1})-\vW_2(t_n).
\]
For the velocity noise define
\begin{equation}\label{eq:Isig}
\Isig^n
:=
P_h\int_{t_n}^{t_{n+1}}
\sigma(\vu^n)\,d\vW_1(s)
\in\bJ_h .
\end{equation}
Since $\vu^n$ is $\mathcal F_{t_n}$-measurable,
\begin{equation}\label{eq:Isig-prop}
\bE\!\left[
\Isig^n\mid\mathcal F_{t_n}
\right]
=0,
\qquad
\bE\!\left[
\norm{\Isig^n}_{L^2}^2
\mid\mathcal F_{t_n}
\right]
\le
\tau\norm{\sigma(\vu^n)}_{\rm HS}^2 .
\end{equation}

We take
\[
\vu^0=P_h\vu_0,
\qquad
\vd^0=R_h\vd_0 .
\]
No normalization of $\vd^0$ is imposed. For $(\vu^n,\vd^n)\in\bJ_h\times\bD_h$, define
\begin{equation}\label{eq:Eh}
\Phi_h(\vd^n)
:=
\frac12\norm{\nabla\vd^n}_{L^2}^2
+
\frac1{\delta^2}\iph{F(\vd^n)}{1},
\qquad
\Ecal_h^n
:=
\frac12\norm{\vu^n}_{L^2}^2+\Phi_h(\vd^n).
\end{equation}

\begin{lemma}[Initial energy]\label{lem:E0}
There is a constant $C$, independent of $h$, such that
\begin{equation}\label{eq:E0}
\Ecal_h^0
\le
C\left(
1+\norm{\vu_0}_{L^2}^2
+\norm{\vd_0}_{\bH^1}^2
+\delta^{-2}\norm{\vd_0}_{\bH^1}^4
\right).
\end{equation}
\end{lemma}

\begin{proof}
The $L^2$-stability of $P_h$ and the $\bH^1$-stability of the Ritz
projection give
\[
\norm{\vu^0}_{L^2}
\le\norm{\vu_0}_{L^2},
\qquad
\norm{\vd^0}_{\bH^1}
\le\norm{\vd_0}_{\bH^1}.
\]
Moreover, by the quartic growth of $F$, \eqref{eq:lumpLq} with $q=4$,
and $\bH^1\hookrightarrow L^4$ in two dimensions,
\[
\iph{F(\vd^0)}{1}
\le
C\bigl(1+\norm{\vd^0}_{L^4}^4\bigr)
\le
C\bigl(1+\norm{\vd_0}_{\bH^1}^4\bigr).
\]
This proves \eqref{eq:E0}.
\end{proof}

\begin{remark}\label{rem:E0}
The estimate above is sufficient for the stability argument; no
$L^q$-stability of the Ritz projection is required.
\end{remark}

\medskip
\noindent
\textbf{Fully discrete scheme.}
Given $(\vu^n,\vd^n)\in\bJ_h\times\bD_h$, find
\[
(\vu^{n+1},p^{n+1},\vd^{n+1})
\in
\bV_h\times\bQ_h\times\bD_h
\]
such that, for every $(v,q,e)\in\bV_h\times\bQ_h\times\bD_h$,
\begin{align}
\ip{\vu^{n+1}-\vu^n}{v}
&+
\tau\widehat b(\vu^n;\vu^{n+1},v)
+
\tau\ip{\nabla\vu^{n+1}}{\nabla v}
-
\tau\ip{p^{n+1}}{\nabla\!\cdot v}
\nonumber\\
&=
-\tau\dual{M_h(\vd^{n+1})}{v}
+
\ip{\Isig^n}{v},
\label{eq:S1}
\\[2mm]
\ip{\nabla\!\cdot\vu^{n+1}}{q}
&=0,
\label{eq:S2}
\\[2mm]
\iph{\vd^{n+1}-\vd^n}{e}
&+
\tau\btil(\vu^{n+1},\vd^{n+1},e)
+
\tau\ip{\nabla\vd^{n+1}}{\nabla e}
\nonumber\\
&\quad
+
\frac{\tau}{\delta^2}
\iph{f_\star(\vd^n,\vd^{n+1})}{e}
=
a_n
\iph{
\dmid\times\vze
}{e},
\label{eq:S3}
\end{align}
where
\[
\dmid:=\frac12(\vd^n+\vd^{n+1}).
\]

Three structural choices are visible in
\eqref{eq:S1}--\eqref{eq:S3}. The director noise is evaluated at the
midpoint, while the dissipative terms are treated by backward Euler; the
latter supplies the numerical dissipation used in the energy estimate.
Moreover, the velocity and director equations are solved monolithically,
which permits the exact cancellation \eqref{eq:BFP}. The compatibility
between the midpoint noise and $f_\star$ is analysed in
\S\ref{sec:structure}.

\subsection{One-step solvability and adaptedness}\label{sub:solv}

\begin{proposition}[Existence of an iterate]\label{prop:solvable}
Fix $n$, $\omega\in\Omega$, $h>0$, $\delta>0$, and $\tau>0$.
Then \eqref{eq:S1}--\eqref{eq:S3} possesses at least one solution
$(\vu^{n+1},p^{n+1},\vd^{n+1})$, with
$\vu^{n+1}\in\bJ_h$.
In particular, one-step solvability requires no relation between
$\tau$, $h$, and $\delta$.
\end{proposition}

\begin{proof}
By the inf--sup condition, it is enough to solve the pressure-free system on
\[
X_h:=\bJ_h\times\bD_h
\]
and recover $p^{n+1}$ afterwards. Equip $X_h$ with
\[
[(v,e),(w,g)]
:=
\ip{v}{w}+\iph{e}{g}.
\]
Let $\mathcal R_u$ and $\mathcal R_d$ denote the residuals of
\eqref{eq:S1} and \eqref{eq:S3}, and define $P:X_h\to X_h$ by
\begin{equation}\label{eq:Pdef}
[P(v,e),(w,g)]
:=
\mathcal R_u(v,e;w)
+
\mathcal R_d\bigl(v,e;(\mathrm{Id}+\Ah)g\bigr).
\end{equation}
Since $\mathrm{Id}+\Ah$ is invertible on $\bD_h$, zeros of $P$ are
precisely the required iterates. The map $P$ is continuous.

Let $z=(v,e)$. By \eqref{eq:BFP},
\[
\tau\dual{M_h(e)}{v}
+
\tau\btil(v,e,\Ah e)=0,
\]
while pointwise divergence freedom gives
\begin{equation}\label{eq:transport-zero}
\btil(v,e,e)
=
\frac12\int_{\bO}v\cdot\nabla\abs e^2\,dx
=0.
\end{equation}
Thus both mixed coupling terms disappear. The remaining coercive terms satisfy
\[
\ip{v-\vu^n}{v}
\ge
\frac12\norm v_{L^2}^2-C_n,
\]
and
\[
\iph{e-\vd^n}{(\mathrm{Id}+\Ah)e}
\ge
\frac12\norm e_{\bH^1}^2-C_n.
\]
Moreover, Lemma \ref{lem:fstarprop}(iv), \eqref{eq:lumpLq}, and Jensen's
inequality yield
\begin{equation}\label{eq:penalty-coercive}
\iph{f_\star(\vd^n,e)}{e}
\ge
c\abs{\bO}^{-1}\norm e_h^4-C_n.
\end{equation}

It remains to control the penalty tested against $\Ah e$. From
\eqref{eq:edge},
\[
\iph{\Ah e}{f_\star(\vd^n,e)}
=
\sum_{i<j}\omega_{ij}
(e_i-e_j)\cdot
\bigl(
f_\star(d_i^n,e_i)-f_\star(d_j^n,e_j)
\bigr).
\]
Split
\[
\begin{aligned}
f_\star(d_i^n,e_i)-f_\star(d_j^n,e_j)
={}&
f_\star(d_i^n,e_i)-f_\star(d_i^n,e_j) +
f_\star(d_i^n,e_j)-f_\star(d_j^n,e_j).
\end{aligned}
\]
The first term is controlled from below by
\eqref{eq:osl}. For the second, Lemma \ref{lem:fstarprop}(ii) and the
fact that $\vd^n$ is fixed give
\[
\abs{
f_\star(d_i^n,e_j)-f_\star(d_j^n,e_j)
}
\le
C_n(1+\abs{e_j}^2)
\abs{d_i^n-d_j^n}.
\]
Since $X_h$ is finite dimensional, norm equivalence then gives
\begin{equation}\label{eq:penalty-Ah}
\iph{\Ah e}{f_\star(\vd^n,e)}
\ge
-C_{n,h}\bigl(1+\norm e_h^3\bigr).
\end{equation}

Finally, the midpoint-noise contribution satisfies $\iph{e\times\vze}{e}=0,$ and, for fixed $h$,
\begin{equation}\label{eq:noise-solvability}
\left|
a_n
\iph{
\tfrac12(\vd^n+e)\times\vze
}{
(\mathrm{Id}+\Ah)e
}
\right|
\le
C_{n,h,\omega}\bigl(1+\norm e_h^2\bigr).
\end{equation}
The velocity noise contributes only a term linear in $v$.

Collecting these estimates,
\begin{equation}\label{eq:brouwer}
\begin{aligned}
[P(z),z]
\ge{}&
c\norm v_{L^2}^2
+c\norm e_{\bH^1}^2
+\tau\norm{\nabla v}_{L^2}^2
+\tau\norm{\Ah e}_h^2
+\frac{c\tau}{\delta^2}\norm e_h^4
\\
&-
C_{n,h,\tau,\delta,\omega}
\bigl(
1+\norm v_{L^2}
+\norm e_h^2
+\norm e_h^3
\bigr).
\end{aligned}
\end{equation}
Hence
\[
[P(z),z]\longrightarrow+\infty
\qquad
\text{as }\norm z_{X_h}\longrightarrow\infty.
\]
The finite-dimensional Brouwer lemma
\cite[Ch.~II, Lemma~1.4]{Temam1979} therefore gives a zero of $P$.
The pressure is then recovered from \eqref{eq:infsup}.
\end{proof}

\begin{remark}[Role of the discrete structure]\label{rem:solv}
Three ingredients are used essentially in the preceding proof:
the energy pairing $(v,(\mathrm{Id}+\Ah)e)$ eliminates the elastic coupling,
the pointwise divergence constraint eliminates
$\btil(v,e,e)$, and weak acuteness permits the edgewise control
\eqref{eq:penalty-Ah}. Thus solvability imposes no mesh-dependent
restriction on $\tau$. Weak acuteness is a sufficient condition for our
argument; we do not claim that the scheme fails on other meshes.
\end{remark}

\begin{proposition}[Adapted discrete solution]\label{prop:unique}
There is a solution $(\vu^n,\vd^n)_{n=0}^{\Nt}$ of
\eqref{eq:S1}--\eqref{eq:S3} that is adapted to the filtration
$(\mathcal F_{t_n})_{n=0}^{\Nt}$.
\end{proposition}

\begin{proof}
Let $Z:=\bJ_h\times\bD_h$ and, for data $z\in Z$ and increments
$\xi=(\theta,a)\in\bJ_h\times\bR$, denote by
$\mathcal S(z,\xi)\subset Z$ the set of solutions of one discrete step.

Proposition \ref{prop:solvable} gives nonempty values. The residual map is
jointly continuous in $(z,\xi,z')$, hence the graph of $\mathcal S$ is
closed. Moreover, the coercivity estimate \eqref{eq:brouwer} holds for all
sufficiently large $\norm{z'}_Z$; consequently every zero lies in a bounded
set. Thus $\mathcal S(z,\xi)$ is compact.

Since the parameter and solution spaces are finite-dimensional Polish spaces,
the measurable-selection theorem
\cite[Thm.~8.1.3]{AubinFrankowska1990}
provides a Borel selector
\[
\mathfrak s:
Z\times(\bJ_h\times\bR)\longrightarrow Z,
\qquad
\mathfrak s(z,\xi)\in\mathcal S(z,\xi).
\]
Set
\[
(\vu^{n+1},\vd^{n+1})
=
\mathfrak s
\bigl(
(\vu^n,\vd^n),(\Isig^n,a_n)
\bigr).
\]
Since $(\Isig^n,a_n)$ is
$\mathcal F_{t_{n+1}}$-measurable, induction yields the required
adaptedness.
\end{proof}

\begin{remark}[Why uniqueness is not required]\label{rem:adapted}
The proof does not require uniqueness of the nonlinear algebraic step.
Measurable selection provides an adapted branch, and every stability,
convergence, and error estimate below is obtained directly from
\eqref{eq:S1}--\eqref{eq:S3}; hence it applies to any solution of the
discrete system. Possible uniqueness under additional small-step conditions
is therefore unnecessary for the analysis.
\end{remark}

\section{A compatibility principle for the noise and the penalty}
\label{sec:structure}

We now identify the structural reason for the midpoint treatment of the
director noise and the $|\vd|^{2}$-secant introduced in
\S\ref{sec:scheme}. The main result is a compatibility principle:
among continuous field-independent two-point discretisations, exact
potential-energy balance and orthogonality to the midpoint geometric noise
uniquely determine the secant. The argument is pointwise and therefore
extends beyond the particular liquid-crystal system considered here.

We begin with three elementary identities. For
$\vze\in\bR^{3}$ let
\[
\Gop:=-[\vze]_{\times},
\qquad
G(\vd)=\Gop\vd=\vd\times\vze .
\]

\begin{lemma}[Algebraic identities]\label{lem:algebraic}
For $\vd\in\bR^{3}$,
\begin{enumerate}[label=\emph{(\roman*)},leftmargin=2.2em,itemsep=1pt]
\item
$\Gop^{\!\top}=-\Gop$ and
$\Gop^{3}=-\abs{\vze}^{2}\Gop$; in particular
\[
G(\vd)\cdot\vd=0.
\]

\item
If $f(\vd)=\phi(\abs{\vd}^{2})\vd$ is radial, then
\[
G(\vd)\cdot f(\vd)=0.
\]

\item
For spatially dependent $\vd$ and $\vze$,
\[
\partial_i(\vd\times\vze)\cdot\partial_i\vd
=
(\vd\times\partial_i\vze)\cdot\partial_i\vd .
\]
Hence this contribution vanishes when $\vze$ is spatially constant.
\end{enumerate}
\end{lemma}

\begin{proof}
The first two assertions follow from skew-symmetry and radiality.
For (iii),
\[
\partial_i(\vd\times\vze)
=
(\partial_i\vd)\times\vze
+
\vd\times\partial_i\vze,
\]
and the first term is orthogonal to $\partial_i\vd$.
\end{proof}

\subsection{The compatibility principle}\label{sub:principle}

The argument is independent of the fluid equation. Consider a director
evolution of the form
\begin{equation}\label{eq:class}
d\vd+
\left[
\mathcal L\vd+
\frac1{\delta^{2}}f(\vd)
\right]dt
=
\sum_{k=1}^{K}
G_k(\vd)\circ d\vW^{k},
\qquad
G_k(\vd)=\vd\times\vze_k,
\end{equation}
where
\begin{equation}\label{eq:radialF}
F(\vd)=\frac12\Ftil(\abs{\vd}^{2}),
\qquad
f(\vd)=F'(\vd)
      =\Ftil'(\abs{\vd}^{2})\vd .
\end{equation}
The spatial operator $\mathcal L$ will play no role below.

\begin{definition}[Evolution with a penalised sphere constraint]
\label{def:class}
We call \eqref{eq:class} an evolution with a penalised sphere constraint if
the potential is radial as in \eqref{eq:radialF} and each noise field is
geometric,
\[
G_k(\vd)\cdot\vd=0 .
\]
The Ginzburg--Landau director equation \eqref{eq:abs-d} is the principal
example considered in this paper.
\end{definition}

Let
\[
\varrho
:=
\frac12(\vd^{n}+\vd^{n+1}),
\qquad
\Delta\vd^{n}
:=
\vd^{n+1}-\vd^{n}.
\]
For a continuous two-point approximation
$f_\star:\bR^{3}\times\bR^{3}\to\bR^{3}$, consider
\begin{align}
\text{\emph{(A) Energy exactness:}}\qquad
&
f_\star(\vd^{n},\vd^{n+1})
\cdot\Delta\vd^{n}
=
F(\vd^{n+1})-F(\vd^{n}),
\label{eq:A}
\\
\text{\emph{(B) Noise orthogonality:}}\qquad
&
(\varrho\times\vze)
\cdot f_\star(\vd^{n},\vd^{n+1})
=0
\quad
\text{for every }\vze\in\bR^{3}.
\label{eq:B}
\end{align}

Condition \emph{(A)} gives an exact discrete chain rule for the potential,
while \emph{(B)} prevents the midpoint noise from entering the penalisation
part of the energy balance.

\begin{theorem}[Compatibility characterisation]
\label{thm:characterisation}
Let $\Ftil\in C^{1}([0,\infty))$ and let
$f_\star:\bR^{3}\times\bR^{3}\to\bR^{3}$ be continuous and independent
of the applied field $\vze$. Then \eqref{eq:A} and \eqref{eq:B} hold for
all $\vd^{n},\vd^{n+1},\vze\in\bR^{3}$ if and only if
\begin{equation}\label{eq:character-sec}
f_\star(\vd^{n},\vd^{n+1})
=
\kappa(\vd^{n},\vd^{n+1})\,\varrho ,
\end{equation}
where
\begin{equation}\label{eq:character-kappa}
\kappa(\xi,\eta)
:=
\begin{cases}
\dfrac{\Ftil(\abs{\eta}^{2})-\Ftil(\abs{\xi}^{2})}
      {\abs{\eta}^{2}-\abs{\xi}^{2}},
&
\abs{\eta}\ne\abs{\xi},
\\[3mm]
\Ftil'(\abs{\xi}^{2}),
&
\abs{\eta}=\abs{\xi}.
\end{cases}
\end{equation}
Thus the admissible two-point nonlinearity is precisely the
$|\vd|^{2}$-secant multiplied by the midpoint.
\end{theorem}

\begin{proof}
Sufficiency follows directly from
\[
\varrho\cdot\Delta\vd^{n}
=
\frac12
\left(
\abs{\vd^{n+1}}^{2}
-
\abs{\vd^{n}}^{2}
\right)
\]
and $(\varrho\times\vze)\cdot\varrho=0$.

For necessity, suppose first that $\varrho\ne0$. Since $\{
\varrho\times\vze:\vze\in\bR^{3}
\} =
\varrho^\perp,$ condition \eqref{eq:B} implies
\[
f_\star(\vd^{n},\vd^{n+1})
\in
(\varrho^\perp)^\perp
=
\operatorname{span}\{\varrho\}.
\]
Hence $f_\star=\kappa\varrho$ for some scalar $\kappa$. Substitution into
\eqref{eq:A} gives
\[
\frac{\kappa}{2}
\left(
\abs{\vd^{n+1}}^{2}
-
\abs{\vd^{n}}^{2}
\right)
=
\frac12
\left[
\Ftil(\abs{\vd^{n+1}}^{2})
-
\Ftil(\abs{\vd^{n}}^{2})
\right],
\]
which yields \eqref{eq:character-kappa} whenever the two squared norms are
different. Since $\Ftil\in C^{1}$ and $f_\star$ is continuous, the same
formula extends uniquely to equal norms and to $\varrho=0$.
\end{proof}

\begin{remark}[Scope of the characterisation]\label{rem:whichP}
Only radiality of the potential and geometricity of the noise are used.
The spatial operator $\mathcal L$, the transport term, and the fluid
equation play no role. The result therefore applies equally to other
penalised sphere-valued evolutions with cross-product Stratonovich noise.

The quantification over every $\vze$ is what forces collinearity with the
midpoint. A single prescribed field leaves additional freedom. More generally,
a sufficiently rich family of physical noise fields yields the same conclusion;
for example, two non-parallel fields determine the result away from the plane
they span, and continuity extends it to the remaining set.
\end{remark}

For later reference we record the two identities separately.

\begin{lemma}[Exactness of the secant]\label{lem:fstar}
Let $f_\star$ be given by \eqref{eq:character-sec}--\eqref{eq:character-kappa}.
Then
\begin{align}
f_\star(\vd^{n},\vd^{n+1})
\cdot\Delta\vd^{n}
&=
F(\vd^{n+1})-F(\vd^{n}),
\tag{A}
\\
(\varrho\times\vze)
\cdot f_\star(\vd^{n},\vd^{n+1})
&=0
\qquad
\text{for every }\vze\in\bR^{3}.
\tag{B}
\end{align}
\end{lemma}

\begin{proof}
Both identities follow immediately from
$f_\star=\kappa\varrho$ and the definition of $\kappa$.
\end{proof}

The second condition has no deterministic counterpart. Without noise,
\eqref{eq:B} disappears and energy exactness alone does not determine a
unique two-point discretisation.

\subsection{Why an energy-exact secant is not enough}

The average-vector-field discretisation illustrates the role of
\eqref{eq:B}. For the quartic Ginzburg--Landau potential, introduce the
midpoint and the corresponding increment
\[
\varrho:=\frac12(\vd^{n}+\vd^{n+1}),
\qquad
\boldsymbol{\eta}^{n}:=\vd^{n+1}-\vd^{n}.
\]

\begin{lemma}[The average-vector-field secant]\label{lem:avfcost}
The average-vector-field discretisation
\[
f_{\rm AVF}
:=
\int_{0}^{1}
f\bigl(\vd^{n}+s\boldsymbol{\eta}^{n}\bigr)\,ds
\]
is energy exact, but
\begin{equation}\label{eq:avf}
f_{\rm AVF}
=
(\abs{\varrho}^{2}-1)\varrho
+
\frac16(\varrho\cdot\boldsymbol{\eta}^{n})\boldsymbol{\eta}^{n}
+
\frac1{12}\abs{\boldsymbol{\eta}^{n}}^{2}\varrho .
\end{equation}
Consequently,
\begin{equation}\label{eq:avf-orth}
(\varrho\times\vze)\cdot f_{\rm AVF}
=
\frac16
(\varrho\cdot\boldsymbol{\eta}^{n})
\bigl((\varrho\times\vze)\cdot\boldsymbol{\eta}^{n}\bigr),
\end{equation}
which is generally non-zero.
\end{lemma}

\begin{proof}
Since $\vd^{n}
=
\varrho-\frac12\boldsymbol{\eta}^{n},$ we may write
\[
\vd^{n}+s\boldsymbol{\eta}^{n}
=
\varrho+t\boldsymbol{\eta}^{n},
\qquad
t=s-\frac12.
\]
Substituting this expression into the definition of $f_{\rm AVF}$ and
integrating over $t\in[-1/2,1/2]$, the odd powers of $t$ vanish, while
\[
\int_{-1/2}^{1/2}t^{2}\,dt=\frac1{12}.
\]
This gives \eqref{eq:avf}. Equation \eqref{eq:avf-orth} then follows from
$(\varrho\times\vze)\cdot\varrho=0$.
\end{proof}

\begin{remark}[Why noise orthogonality matters]\label{rem:whyB}
The transverse component $\frac16
(\varrho\cdot\boldsymbol{\eta}^{n})\boldsymbol{\eta}^{n}$ in \eqref{eq:avf} is precisely the term that violates the noise-orthogonality
condition \eqref{eq:B}. When inserted into the penalised energy balance, it
produces a stochastic residual carrying the factor $\delta^{-2}$. In contrast,
for the $|\vd|^{2}$-secant this contribution vanishes identically, rather than
merely tending to zero as $\tau\to0$.
\end{remark}

\subsection{Geometry and strong order of the midpoint noise}

Consider first the isolated director-noise equation at one spatial point,
\[
d\vd
=
\Gop\vd\circ d\vW_{2},
\qquad
\Gop=-[\vze]_{\times}.
\]
The midpoint step is
\begin{equation}\label{eq:noise-midstep}
\vd^{n+1}-\vd^{n}
=
a_n\Gop\varrho,
\qquad
a_n:=\vW_{2}(t_{n+1})-\vW_{2}(t_n).
\end{equation}

\begin{proposition}[Cayley structure and strong order one]
\label{prop:cayley}
Let $I$ denote the identity on $\bR^{3}$. Then:

\begin{enumerate}[label=\emph{(\roman*)},leftmargin=2.2em,itemsep=2pt]

\item
For every $a\in\bR$, $I-\frac a2\Gop$ is invertible and
\begin{equation}\label{eq:cayley}
\vd^{n+1}
=
R(a_n)\vd^{n},
\qquad
R(a)
=
\left(I-\frac a2\Gop\right)^{-1}
\left(I+\frac a2\Gop\right).
\end{equation}

\item
$R(a)$ is orthogonal: $R(a)^{\!\top}R(a)=I\,.$ Consequently,
\[
\abs{\vd^{n+1}}
=
\abs{\vd^{n}}
\]
pathwise for every time step. With the lumped mass matrix this identity
holds independently at every finite element node.

\item
The Cayley transform satisfies
\begin{equation}\label{eq:cayley-exp}
R(a)-e^{a\Gop}
=
\frac1{12}a^{3}\Gop^{3}
+
O(a^{4}).
\end{equation}
Since $\norm{a_n}_{L^{p}(\Omega)}
= O(\tau^{1/2}),$ the one-step defect is $O(\tau^{3/2})$ in every finite $L^{p}(\Omega)$.

\item
Its It\^o expansion is
\begin{equation}\label{eq:cayley-ito}
R(a_n)
=
I+a_n\Gop
+\frac{\tau}{2}\Gop^{2}
+\frac{a_n^{2}-\tau}{2}\Gop^{2}
+O(a_n^{3}).
\end{equation}
Thus the midpoint rule reproduces the scalar Milstein correction and has
strong order one, without simulation of a L\'evy area.
\end{enumerate}
\end{proposition}

\begin{proof}
Since $\Gop$ is skew, its spectrum is purely imaginary and
$I-\frac a2\Gop$ is invertible. Orthogonality of \eqref{eq:cayley} is the
standard Cayley-transform identity; equivalently, testing
\eqref{eq:noise-midstep} with $\varrho$ gives
\[
\abs{\vd^{n+1}}^{2}-\abs{\vd^{n}}^{2}
=
2a_n(\Gop\varrho)\cdot\varrho
=0.
\]
Expanding the Cayley transform and the exponential gives
\eqref{eq:cayley-exp}. Finally,
\[
\int_{t_n}^{t_{n+1}}\!\int_{t_n}^{s}
d\vW_{2}(r)\,d\vW_{2}(s)
=
\frac12(a_n^{2}-\tau),
\]
which identifies the last two terms in \eqref{eq:cayley-ito} with the
It\^o--Stratonovich and Milstein corrections; see
\cite[Ch.~10]{KloedenPlaten1992}.
\end{proof}

\subsection{Several director noises}\label{sub:multi}

The geometric property survives several independent driving Wiener processes,
whereas the strong order depends on commutativity.

\begin{proposition}[Multi-dimensional director noise]\label{prop:multi}
Consider
\[
d\vd
=
\sum_{k=1}^{K}
\Gop_k\vd\circ d\vW^{k},
\qquad
\Gop_k=-[\vze_k]_{\times}.
\]
Then:

\begin{enumerate}[label=\emph{(\roman*)},leftmargin=2.2em,itemsep=2pt]

\item
we have
\[
[\Gop_k,\Gop_l]
=
[\vze_k\times\vze_l]_{\times},
\]
so the operators commute if and only if the corresponding fields are
pairwise parallel.

\item
The compatibility principle of Theorem \ref{thm:characterisation} is
unchanged.

\item
The midpoint step
\[
\vd^{n+1}-\vd^{n}
=
\sum_{k=1}^{K}
a_n^{k}\Gop_k\varrho
\]
is the Cayley transform of the skew matrix
$\sum_k a_n^{k}\Gop_k$ and is therefore an exact isometry, irrespective
of commutativity.

\item
If the fields are pairwise parallel, the noise reduces to a scalar Brownian
motion multiplying a fixed skew operator, and the midpoint rule has strong
order one. For genuinely non-commuting fields, order-one approximation
requires the corresponding L\'evy-area terms; without them the midpoint
method has, in general, strong order $1/2$.
\end{enumerate}
\end{proposition}

\begin{proof}
The commutator identity follows from
\[
[[\xi]_{\times},[\eta]_{\times}]
=
[\xi\times\eta]_{\times}.
\]
A sum of skew matrices is skew, which proves the isometry statement.
For commuting fields the problem reduces to the scalar situation of
Proposition \ref{prop:cayley}. In the non-commuting case, the stochastic
Taylor expansion contains commutators multiplied by L\'evy areas, which
are absent from the midpoint update.
\end{proof}

\begin{remark}\label{rem:multi}
Thus geometry and temporal order should be distinguished: the midpoint
method preserves the quadratic constraint for any number of geometric
noises, while strong order one without L\'evy-area simulation is recovered
only in the commuting case.
\end{remark}

\subsection{The cost of treating the noise explicitly}

We finally compare the midpoint rule with the Euler--Maruyama treatment of
the equivalent It\^o equation. The mechanism is already visible in the
isolated noise subproblem, so no deterministic terms are needed.

Let $x_h\in\bD_h$ be $\mathcal F_{t_n}$-measurable and let $\vze$ be
spatially constant. The midpoint update is
\[
y_h^{\rm mid}-x_h
=
a_nG\!\left(\frac{x_h+y_h^{\rm mid}}2\right),
\]
whereas the explicit It\^o-corrected update is
\begin{equation}\label{eq:explicit-noise-step}
y_h^{\exp}
=
x_h
+
a_nG(x_h)
+
\frac{\tau}{2}G^{2}(x_h).
\end{equation}
Both are paired with the same secant $f_\star$.

Define the penalty contributions
\begin{align}
\mathcal N_{\rm mid}^{n}
&:=
\frac{a_n}{\delta^{2}}
\iph{
G\!\left(\frac{x_h+y_h^{\rm mid}}2\right)
}{
f_\star(x_h,y_h^{\rm mid})
},
\label{eq:NGLmid}
\\
\mathcal N_{\exp}^{n}
&:=
\frac1{\delta^{2}}
\iph{
a_nG(x_h)+\frac{\tau}{2}G^{2}(x_h)
}{
f_\star(x_h,y_h^{\exp})
}.
\label{eq:NGL}
\end{align}

\begin{proposition}[Explicit noise: bias and fluctuation]
\label{prop:negative}
Assume that the iterates of \eqref{eq:explicit-noise-step} possess moments
of sufficiently high order uniformly on $[0,T]$. Then:

\begin{enumerate}[label=\emph{(\roman*)},leftmargin=2.2em,itemsep=2pt]

\item We have $\mathcal N_{\rm mid}^{n}=0$ identically, pathwise and for every $\tau$ and $\delta$.

\item
Let $\kappa_h
:= \kappa(x_h,y_h^{\exp})$ be understood nodewise. Then
\begin{equation}\label{eq:Nexp-exact}
\mathcal N_{\exp}^{n}
=
\frac{a_n^{2}-\tau}{2\delta^{2}}
\iph{\kappa_h G(x_h)}{G(x_h)}
+
\frac{\tau^{2}}{8\delta^{2}}
\iph{\kappa_h G^{2}(x_h)}{G^{2}(x_h)} .
\end{equation}

\item
Over $N=T/\tau$ steps,
\begin{equation}\label{eq:Nexp-scales}
\left|
\bE\sum_{n=0}^{N-1}\mathcal N_{\exp}^{n}
\right|
\le
C\delta^{-2}\tau,
\qquad
\left(
\operatorname{Var}
\sum_{n=0}^{N-1}\mathcal N_{\exp}^{n}
\right)^{1/2}
\le
C\delta^{-2}\tau^{1/2},
\end{equation}
where $C$ is independent of $\tau$ and $\delta$. Consequently,
the second-moment estimate is $\delta$-uniform only under the scaling $\tau=O(\delta^{4}).$
\end{enumerate}
\end{proposition}

\begin{proof}
Part (i) is Lemma \ref{lem:fstar}(B). For the explicit update, put
\[
q_h
:=
a_nG(x_h)+\frac{\tau}{2}G^{2}(x_h),
\qquad
\varrho_h^{\exp}
=
x_h+\frac12q_h.
\]
Since $\Gop$ is skew,
\[
G(x_h)\cdot x_h=0,
\qquad
G(x_h)\cdot G^{2}(x_h)=0,
\qquad
G^{2}(x_h)\cdot x_h
=
-\abs{G(x_h)}^{2}.
\]
Hence, nodewise,
\[
q_h\cdot\varrho_h^{\exp}
=
\frac12(a_n^{2}-\tau)\abs{G(x_h)}^{2}
+
\frac{\tau^{2}}8\abs{G^{2}(x_h)}^{2}.
\]
Since
$f_\star(x_h,y_h^{\exp})=\kappa_h\varrho_h^{\exp}$,
this proves \eqref{eq:Nexp-exact}.

For the quartic potential,
\[
\kappa_h
=
\ftil(\abs{x_h}^{2})
+
\frac12(a_n^{2}-\tau)\abs{G(x_h)}^{2}
+
\frac{\tau^{2}}8\abs{G^{2}(x_h)}^{2}
\]
nodewise. Since
\[
\bE[a_n^{2}-\tau]=0,
\qquad
\bE[(a_n^{2}-\tau)^{2}]
=
2\tau^{2},
\]
the conditional mean of one step is
$O(\delta^{-2}\tau^{2})$ and its martingale part has conditional variance
$O(\delta^{-4}\tau^{2})$. Summing $T/\tau$ steps gives
\eqref{eq:Nexp-scales}.
\end{proof}

\begin{remark}[What is not preserved]\label{rem:notpreserved}
Proposition \ref{prop:cayley} concerns the isolated noise problem. In the
full scheme \eqref{eq:S3}, diffusion, transport, and penalisation are present,
so in general
\[
\abs{\vd^{n+1}}
\ne
\abs{\vd^{n}}.
\]
What the midpoint rule guarantees is more precise: the geometric noise itself
does not contribute to the change in the director length or to the penalisation
part of the energy. Any deviation from the sphere is therefore produced by
the deterministic part of the penalised dynamics, as intended.
\end{remark}

The preservation of quadratic invariants by the implicit midpoint rule for
Stratonovich systems is classical
\cite{MilsteinRepinTretyakov2002,HairerLubichWanner2006}. The specific feature
used here is its compatibility with the Ginzburg--Landau penalisation:
the midpoint geometry and the $|\vd|^{2}$-secant together remove the
potentially singular $\delta^{-2}$ stochastic contribution exactly.

\section{Discrete energy law and stability}\label{sec:stability}

The compatibility principle of Theorem \ref{thm:characterisation} yields an
exact discrete energy identity. We first derive this identity, then isolate the
only residual produced by the transport--penalty interaction, and finally
obtain localized stability estimates with no coupling between $\tau$ and $h$.

Throughout this section we write
\[
\mathfrak d\vu^{n}:=\vu^{n+1}-\vu^{n},
\qquad
\mathfrak d\vd^{n}:=\vd^{n+1}-\vd^{n},
\qquad
\dmid:=\frac12(\vd^{n+1}+\vd^{n}),
\]
and
\[
w^{n}:=\abs{\dmid}^{2}.
\]
The discrete energy $\Ecal_h^n$ is defined in \eqref{eq:Eh}. We also use the
nonnegative shifted energy
\begin{equation}\label{eq:shifted-energy}
\widehat{\Ecal}_h^n
:=
\Ecal_h^n+C_F\abs{\bO}\delta^{-2}\ge0,
\end{equation}
where $C_F$ is the constant of Lemma \ref{lem:coerc}. Unless otherwise
specified, $C$ denotes a constant independent of $h$, $\tau$, $\delta$ and
$n$.

\subsection{The exact discrete energy balance}

\begin{theorem}[Discrete energy identity]\label{thm:energy}
Let $(\vu^{n+1},p^{n+1},\vd^{n+1})$ satisfy
\eqref{eq:S1}--\eqref{eq:S3}. Then
\begin{multline}\label{eq:energylaw}
\Ecal_h^{n+1}-\Ecal_h^n
+\frac12\norm{\mathfrak d\vu^n}_{L^2}^{2}
+\frac12\norm{\nabla\mathfrak d\vd^n}_{L^2}^{2}
+\tau\norm{\nabla\vu^{n+1}}_{L^2}^{2}
\\
+\tau
\left\lVert
\Ah\vd^{n+1}+\frac1{\delta^2}f_\star
\right\rVert_h^{2}
=
\mathcal M^n+\mathcal I_1^n+\mathcal I_2^n
+\mathcal N^n+\mathcal R^n ,
\end{multline}
where
\begin{align}
\mathcal M^n
&:=\ip{\Isig^n}{\vu^n},
&
\mathcal I_1^n
&:=\ip{\Isig^n}{\mathfrak d\vu^n},
\label{eq:energyterms1}
\\
\mathcal I_2^n
&:=
\frac{a_n}{2}
\ip{\nabla I_h(\dmid\times\vze)}
    {\nabla\mathfrak d\vd^n},
&
\mathcal N^n
&:=a_n\Psi_h(\dmid),
\label{eq:Nn}
\\
\mathcal R^n
&:=
-\frac{\tau}{\delta^2}
\btil(\vu^{n+1},\vd^{n+1},f_\star),
\label{eq:Rdef}
\end{align}
and
\begin{equation}\label{eq:Psih}
\Psi_h(v)
:=
\ip{\nabla I_h(v\times\vze)}{\nabla v},
\qquad v\in\bD_h.
\end{equation}
Moreover,
\begin{equation}\label{eq:Psih-props}
\Psi_h\equiv0
\quad\text{if $\vze$ is spatially constant},
\qquad
\abs{\Psi_h(v)}
\le
C\norm{\vze}_{W^{1,\infty}}\norm{v}_{H^1}^{2}.
\end{equation}

In particular, the stochastic contribution from the Ginzburg--Landau
penalty vanishes identically:
\begin{equation}\label{eq:penalty-noise-zero}
\frac{a_n}{\delta^2}
\iph{\dmid\times\vze}{f_\star}=0 .
\end{equation}
Thus no negative power of $\delta$ enters the stochastic part of the
discrete energy law.
\end{theorem}

\begin{proof}
Take $v=\vu^{n+1}$ in \eqref{eq:S1},
$q=p^{n+1}$ in \eqref{eq:S2}, and $e=\Ah\vd^{n+1}+\delta^{-2}f_\star$ in \eqref{eq:S3}, and add the resulting identities.

For the velocity,
\[
\ip{\mathfrak d\vu^n}{\vu^{n+1}}
=
\frac12
\left(
\norm{\vu^{n+1}}_{L^2}^{2}
-\norm{\vu^n}_{L^2}^{2}
+\norm{\mathfrak d\vu^n}_{L^2}^{2}
\right).
\]
The convection term vanishes by skew-symmetry, the pressure term by
\eqref{eq:S2}, and
\[
\ip{\Isig^n}{\vu^{n+1}}
=
\mathcal M^n+\mathcal I_1^n.
\]

For the director,
\[
\iph{\mathfrak d\vd^n}{\Ah\vd^{n+1}}
=
\frac12
\left(
\norm{\nabla\vd^{n+1}}_{L^2}^{2}
-\norm{\nabla\vd^n}_{L^2}^{2}
+\norm{\nabla\mathfrak d\vd^n}_{L^2}^{2}
\right),
\]
while the exact secant identity of Lemma \ref{lem:fstar} gives
\[
\frac1{\delta^2}
\iph{\mathfrak d\vd^n}{f_\star}
=
\frac1{\delta^2}
\iph{F(\vd^{n+1})-F(\vd^n)}{1}.
\]
The diffusion and penalty terms therefore combine into $\tau
\left\lVert
\Ah\vd^{n+1}+\delta^{-2}f_\star
\right\rVert_h^2.$

The elastic contribution cancels with
$\tau\btil(\vu^{n+1},\vd^{n+1},\Ah\vd^{n+1})$
by \eqref{eq:BFP}; the remaining transport contribution is
$-\mathcal R^n$.

Finally, Lemma \ref{lem:fstar}\emph{(B)} gives
\eqref{eq:penalty-noise-zero}. Hence the director noise contributes only
\[
a_n
\iph{\dmid\times\vze}{\Ah\vd^{n+1}}
=
a_n
\ip{\nabla I_h(\dmid\times\vze)}
    {\nabla\vd^{n+1}} .
\]
Since
\[
\nabla\vd^{n+1}
=
\nabla\dmid+\frac12\nabla\mathfrak d\vd^n,
\]
this is $\mathcal N^n+\mathcal I_2^n$.

If $\vze$ is constant, then
$I_h(v\times\vze)=v\times\vze$ and
$\Psi_h(v)=0$ by Lemma \ref{lem:algebraic}\emph{(iii)}.
The general bound follows from Lemma \ref{lem:disc}\emph{(vi)}.
\end{proof}

\begin{remark}\label{rem:interp}
For nonconstant $\vze$, the product $v_h\times\vze$ need not belong to
$\bD_h$, which is why the exact energy identity naturally contains
$I_h(v_h\times\vze)$. The form \eqref{eq:Psih} is therefore the appropriate
discrete analogue of the corresponding continuum term.
\end{remark}

\subsection{Increment and transport--penalty residual}

We first record a bound for the director increment.

\begin{lemma}[Director increment]\label{lem:incr}
Set
\[
\Xi^n
:=
\norm{\vu^{n+1}}_{L^4}
\norm{\nabla\vd^{n+1}}_{L^4}
+
\delta^{-2}\norm{f_\star}_h .
\]
Then
\begin{equation}\label{eq:incr}
\norm{\mathfrak d\vd^n}_h^2
\le
C\Bigl[
\tau
\norm{\nabla\vd^{n+1}}_{L^2}
\norm{\nabla\mathfrak d\vd^n}_{L^2}
+
\tau^2(\Xi^n)^2
+
a_n^2\norm{\vze}_{L^\infty}^{2}
\norm{\dmid}_h^2
\Bigr].
\end{equation}
\end{lemma}

\begin{proof}
Test \eqref{eq:S3} with $e=\mathfrak d\vd^n$. The diffusion term is bounded
by
\[
\tau
\norm{\nabla\vd^{n+1}}_{L^2}
\norm{\nabla\mathfrak d\vd^n}_{L^2},
\]
whereas the transport and penalty terms are bounded by $\tau\Xi^n\norm{\mathfrak d\vd^n}_h.$ The noise term is bounded by
\[
\abs{a_n}\norm{\vze}_{L^\infty}
\norm{\dmid}_h\norm{\mathfrak d\vd^n}_h.
\]
Young's inequality yields \eqref{eq:incr}.
\end{proof}

The only deterministic remainder in \eqref{eq:energylaw} comes from the
fact that the continuum cancellation
$\btil(\vu,\vd,f(\vd))=0$ is not reproduced exactly by the two-point
nonlinearity.

\begin{lemma}[Structure of the residual]\label{lem:R}
For the quartic potential,
\begin{equation}\label{eq:Rsplit}
\mathcal R^n
=
-\frac{\tau}{\delta^2}
\bigl(T_1^n+T_3^n\bigr),
\end{equation}
where
\begin{align}
T_1^n
&:=
\frac12
\int_{\bO}
\bigl(\kappa-\ftil(w^n)\bigr)
\vu^{n+1}\cdot\nabla w^n\,dx,
\label{eq:T1}
\\
T_3^n
&:=
\frac12
\int_{\bO}
\kappa\,
(\vu^{n+1}\cdot\nabla)\mathfrak d\vd^n
\cdot\dmid\,dx .
\label{eq:T3}
\end{align}
Moreover,
\begin{equation}\label{eq:kappa-mid}
\kappa-\ftil(w^n)
=
\frac14\abs{\mathfrak d\vd^n}^{2},
\end{equation}
and
\begin{equation}\label{eq:T3f}
T_3^n
=
\frac12
\int_{\bO}
(\vu^{n+1}\cdot\nabla)\mathfrak d\vd^n
\cdot f_\star\,dx .
\end{equation}

Consequently, for a fixed polynomial $\Theta^n$ in
$\norm{\vd^n}_{H^1}$ and $\norm{\vd^{n+1}}_{H^1}$,
\begin{align}
\abs{T_1^n}
&\le
C\Theta^n
\norm{\vu^{n+1}}_{L^8}
\norm{\mathfrak d\vd^n}_{L^2}^{1/2}
\norm{\mathfrak d\vd^n}_{H^1}^{3/2},
\label{eq:T1bound}
\\
\abs{T_3^n}
&\le
C\norm{\vu^{n+1}}_{L^4}
\norm{\nabla\mathfrak d\vd^n}_{L^2}
\norm{f_\star}_{L^4}.
\label{eq:T3bound}
\end{align}
No inverse inequality is used in these estimates.
\end{lemma}

\begin{proof}
Since
\[
\vd^{n+1}
=
\dmid+\frac12\mathfrak d\vd^n,
\qquad
f_\star=\kappa\dmid,
\]
we obtain
\[
\btil(\vu^{n+1},\vd^{n+1},f_\star)
=
\frac12
\int_{\bO}
\kappa\,\vu^{n+1}\cdot\nabla w^n\,dx
+
T_3^n .
\]
Add and subtract $\ftil(w^n)$. Since
$\nabla\!\cdot\vu^{n+1}=0$ pointwise,
\[
\frac12
\int_{\bO}
\ftil(w^n)\,
\vu^{n+1}\cdot\nabla w^n\,dx
=
\int_{\bO}
\vu^{n+1}\cdot\nabla
\left[\frac12\Ftil(w^n)\right]dx
=0,
\]
which proves \eqref{eq:Rsplit}.

For the quartic,
\[
\kappa
=
\frac12
\bigl(
\abs{\vd^{n+1}}^2+\abs{\vd^n}^2
\bigr)-1
=
w^n+\frac14\abs{\mathfrak d\vd^n}^2-1,
\]
and hence \eqref{eq:kappa-mid}. Equation \eqref{eq:T3f} follows immediately
from $f_\star=\kappa\dmid$.

For $T_1^n$, use
$\abs{\nabla w^n}\le2\abs{\dmid}\abs{\nabla\dmid}$ and H\"older with
exponents $(8,8,8,8,2)$, followed by the two-dimensional interpolation
estimate
\[
\norm{z}_{L^8}
\le
C\norm{z}_{L^2}^{1/4}\norm{z}_{H^1}^{3/4}.
\]
This gives \eqref{eq:T1bound}. Estimate \eqref{eq:T3bound} follows directly
from H\"older with exponents $(4,2,4)$.
\end{proof}

\begin{remark}[The missing continuum cancellation]\label{rem:fourth}
The first three cancellations in \eqref{eq:cont-cancel} survive exactly at
the discrete level: skew-symmetry of the convection form, pointwise
divergence freedom, and the elastic--transport cancellation
\eqref{eq:BFP}. The fourth,
\[
\btil(\vu,\vd,f(\vd))=0,
\]
does not. Its defect is precisely $\mathcal R^n$. The two terms in
\eqref{eq:Rsplit} have different origins: $T_1^n$ measures the difference
between the temporal secant $\kappa$ and the spatial chain rule, while
$T_3^n$ results from evaluating transport at $\vd^{n+1}$ whereas
$f_\star$ is aligned with the midpoint.
\end{remark}

\begin{remark}\label{rem:nomesh}
Estimates \eqref{eq:T1bound}--\eqref{eq:T3bound} use only
two-dimensional Sobolev interpolation and no inverse estimate. Therefore the
stability argument below requires no condition coupling $\tau$ to $h$.
\end{remark}

\subsection{Localized stability}

Because the residual contains polynomial functions of the discrete solution,
we localize the estimates. For $K>0$ define
\begin{equation}\label{eq:stopindex}
n_K
:=
\min
\left\{
n\le\Nt:
\widehat{\Ecal}_h^n>K
\right\}
\wedge\Nt .
\end{equation}

The following estimate summarizes the higher-order terms generated by
$\mathcal I_2^n$, $\mathcal N^n$ and $\mathcal R^n$.

\begin{proposition}[Localized energy estimate]\label{prop:pathwise}
There exists a polynomial $P$, independent of $h$, $\tau$, $\delta$ and
$K$, such that for every $K>0$ and every $\ell\le n_K$,
\begin{align}
&\bE\Bigl[
\max_{m\le\ell}\widehat{\Ecal}_h^m
\Bigr]
+
c\,\bE\sum_{n<\ell}
\norm{\nabla\mathfrak d\vd^n}_{L^2}^{2} +
c\tau\,\bE\sum_{n<\ell}
\left(
\norm{\nabla\vu^{n+1}}_{L^2}^{2}
+
\left\lVert
\Ah\vd^{n+1}+\frac{f_\star}{\delta^2}
\right\rVert_h^{2}
\right) \nonumber
\\
&\le
C\bigl(1+\widehat{\Ecal}_h^0\bigr)
+
C\tau\sum_{n<\ell}
\bE\Bigl[
1+\widehat{\Ecal}_h^n
\Bigr]
+
C\,T\,\tau\,
\Gamma(K,\delta),
\label{eq:pathwise}
\end{align}
where one may take
\begin{equation}\label{eq:Gamma}
\Gamma(K,\delta)
=
(1+\delta^{-8})P(K).
\end{equation}
The constants $c,C>0$ are independent of $h,\tau,K$ and $\delta$.
\end{proposition}

\begin{proof}
Sum \eqref{eq:energylaw} up to $\ell\le n_K$. For $\mathcal I_1^n$, Young's inequality gives
\[
\mathcal I_1^n
\le
\frac12\norm{\mathfrak d\vu^n}_{L^2}^{2}
+
\frac12\norm{\Isig^n}_{L^2}^{2},
\]
and
\[
\bE\!\left[
\norm{\Isig^n}_{L^2}^{2}
\mid\mathcal F_{t_n}
\right]
\le
C\tau
\bigl(1+\norm{\vu^n}_{L^2}^{2}\bigr).
\]

For $\mathcal I_2^n$,
\[
\abs{\mathcal I_2^n}
\le
\frac1{16}
\norm{\nabla\mathfrak d\vd^n}_{L^2}^{2}
+
C a_n^2\norm{\dmid}_{H^1}^{2}.
\]
Writing $\dmid=\vd^n+\frac12\mathfrak d\vd^n$ and using
Lemma \ref{lem:incr}, the terms containing the increment are either absorbed
into the dissipation or contribute, after expectation and stopping, at most
\[
C\tau\bigl(1+\widehat{\Ecal}_h^n\bigr)
+
C\tau^2P(K).
\]

Next, since $\Psi_h$ is quadratic,
\[
\Psi_h(\dmid)
=
\Psi_h(\vd^n)
+
B_h(\vd^n,\mathfrak d\vd^n)
+
\frac14\Psi_h(\mathfrak d\vd^n),
\]
where
\[
B_h(v,z)
:=
\frac12
\left[
\ip{\nabla I_h(v\times\vze)}{\nabla z}
+
\ip{\nabla I_h(z\times\vze)}{\nabla v}
\right].
\]
Hence, $a_n\Psi_h(\vd^n)$ is a martingale increment, whereas the remaining terms are treated by
Lemma \ref{lem:incr}, Young's inequality and
$\bE a_n^{2m}=C_m\tau^m$. Up to the stopping index they contribute
\[
\varepsilon\,
\bE\sum_{n<\ell}
\norm{\nabla\mathfrak d\vd^n}_{L^2}^{2}
+
C\tau\sum_{n<\ell}
\bE\bigl(1+\widehat{\Ecal}_h^n\bigr)
+
CT\tau P(K).
\]

Finally, Lemma \ref{lem:R} and Young's inequality give
\begin{align*}
\frac{\tau}{\delta^2}\abs{T_3^n}
&\le
\frac1{16}
\norm{\nabla\mathfrak d\vd^n}_{L^2}^{2}
+
\frac{\tau}{4}
\norm{\nabla\vu^{n+1}}_{L^2}^{2}
+
C\tau^3\delta^{-8}P(K),
\\
\frac{\tau}{\delta^2}\abs{T_1^n}
&\le
\frac1{16}
\norm{\mathfrak d\vd^n}_{H^1}^{2}
+
C\tau^4\delta^{-8}P(K).
\end{align*}
The $L^2$ part of the last increment is controlled by
Lemma \ref{lem:incr}. Since $\tau\le1$, summing these higher-order terms
produces at most $CT\tau\,\delta^{-8}P(K).$

The martingales
\[
\sum_{n<\ell}\mathcal M^n,
\qquad
\sum_{n<\ell}a_n\Psi_h(\vd^n)
\]
are estimated by the discrete Burkholder--Davis--Gundy inequality and Young's
inequality. Their contribution is bounded by
\[
\frac14
\bE\max_{m\le\ell}\widehat{\Ecal}_h^m
+
C\tau
\sum_{n<\ell}
\bE\bigl(1+\widehat{\Ecal}_h^n\bigr).
\]
Absorbing the small terms yields \eqref{eq:pathwise}.
\end{proof}

\begin{theorem}[Localized stability]\label{thm:stability}
Let Assumptions \ref{ass:noise}--\ref{ass:F} and \ref{ass:disc} hold and let
$(\vu_0,\vd_0)\in H\times\bH^1$. For every $K>0$ there exists
$\tau_1(K,\delta)>0$, independent of $h$, such that for
$0<\tau\le\tau_1(K,\delta)$,
\begin{equation}\label{eq:stab}
\bE\Bigl[
\max_{0\le n\le n_K}\widehat{\Ecal}_h^n
\Bigr]
+
\bE\sum_{n<n_K}
\norm{\nabla\mathfrak d\vd^n}_{L^2}^{2}
+
\tau\bE\sum_{n<n_K}
\left(
\norm{\nabla\vu^{n+1}}_{L^2}^{2}
+
\left\lVert
\Ah\vd^{n+1}+\frac{f_\star}{\delta^2}
\right\rVert_h^{2}
\right)
\le C_0,
\end{equation}
where
\begin{equation}\label{eq:C0}
C_0
\le
C\bigl(1+\widehat{\Ecal}_h^0\bigr)e^{CT}
\end{equation}
and $C$ is independent of $h$, $\tau$, $K$ and $\delta$.

In particular,
\begin{equation}\label{eq:boundedprob}
\bP\left(
\max_{0\le n\le\Nt}\widehat{\Ecal}_h^n>K
\right)
\le
\frac{C_0}{K},
\end{equation}
and hence the discrete energies are bounded in probability uniformly in
$h$ and $\tau$.
\end{theorem}

\begin{proof}
Apply Proposition \ref{prop:pathwise} with $\ell=n_K$. Choose
\[
\tau_1(K,\delta)
:=
\min\left\{
1,\,
\frac{1}{2CT\,\Gamma(K,\delta)}
\right\}.
\]
Then the last term of \eqref{eq:pathwise} is uniformly bounded. The remaining
linear term is treated by the discrete Gr\"onwall lemma, giving
\eqref{eq:stab}--\eqref{eq:C0}. Importantly, the polynomial dependence on
$K$ and the negative powers of $\delta$ enter the time-step threshold, not
the Gr\"onwall exponent.

Since $\widehat{\Ecal}_h^n\ge0$, Markov's inequality gives
\eqref{eq:boundedprob}.
\end{proof}

\begin{remark}[Uniform in $\delta$]
\label{rem:deltaunif}
Three statements should be distinguished. First, the \emph{stochastic part} of the exact energy identity is free of
negative powers of $\delta$. The potentially dangerous term
\[
\frac{a_n}{\delta^2}
\iph{\dmid\times\vze}{f_\star}
\]
vanishes identically and pathwise by the compatibility principle.

Second, the complete energy identity is not $\delta$-uniform. The
dissipation contains
\[
\tau\bigl\lVert\Ah\vd^{n+1}+\delta^{-2}f_\star\bigr\rVert_h^2,
\]
and the deterministic residual $\mathcal R^n$ carries $\delta^{-2}$.

Third, the stability constant cannot in general be uniform as
$\delta\to0$, because already
\[
\widehat{\Ecal}_h^0
=
\Ecal_h^0+C_F\abs{\bO}\delta^{-2}
\]
may grow like $\delta^{-2}$ for non-well-prepared initial data. What the
argument shows is that no additional negative power of $\delta$ enters the
Gr\"onwall exponent: the stronger $\delta$-dependence occurs only in the
localization threshold $\tau_1(K,\delta)$. Whether that dependence can be
improved under additional structure or well-prepared data is left open.
\end{remark}

\section{Convergence to a martingale solution}\label{sec:convergence}

We now use the stability estimates of Section \ref{sec:stability} to pass to
the limit as $h,\tau\to0$, with the penalty parameter $\delta>0$ fixed.
The main ingredients are compactness in space, fractional regularity in time,
and a careful identification of the midpoint Stratonovich noise.

\subsection{Interpolants and uniform estimates}

For $t\in(t_n,t_{n+1}]$ define
\[
\bar\vu_{h\tau}(t):=\vu^{n+1},
\qquad
\underline\vu_{h\tau}(t):=\vu^n,
\qquad
\hat\vu_{h\tau}(t)
:=
\vu^n+\frac{t-t_n}{\tau}\mathfrak d\vu^n,
\]
and analogously
$\bar\vd_{h\tau}$, $\underline\vd_{h\tau}$ and
$\hat\vd_{h\tau}$, where
\[
\mathfrak d\vu^n:=\vu^{n+1}-\vu^n,
\qquad
\mathfrak d\vd^n:=\vd^{n+1}-\vd^n.
\]

For $K>0$ let
\begin{equation}\label{eq:OmegaK}
\Omega_K:=\{n_K=\Nt\}.
\end{equation}
By Theorem \ref{thm:stability},
\begin{equation}\label{eq:OmegaKprob}
\bP(\Omega_K^c)
\le
\frac{C_0}{K},
\qquad
\tau\le\tau_1(K,\delta),
\end{equation}
where $C_0$ is independent of $K$, $h$ and $\tau$.

\begin{remark}[Order of the localisation limits]\label{rem:order}
Given $\varepsilon>0$, we first choose $K$ such that
$C_0/K<\varepsilon$, and then let $h,\tau\to0$ with $K$ fixed. Since
$\tau\to0$, the restriction
$\tau\le\tau_1(K,\delta)$ is eventually satisfied. Thus the dependence of
$\tau_1$ on $K$ does not obstruct removal of the localisation.
\end{remark}

We shall also retain the velocity-increment term from the proof of
Theorem \ref{thm:stability}; namely,
\begin{equation}\label{eq:u-increments}
\bE\!\left[
\mathbf 1_{\Omega_K}
\sum_{n=0}^{\Nt-1}
\norm{\mathfrak d\vu^n}_{L^2}^{2}
\right]
\le C(K).
\end{equation}

\begin{lemma}[Uniform discrete bounds]\label{lem:unif}
For every fixed $K>0$ there is $C(K,\delta)$, independent of
$h$ and $\tau$, such that
\begin{align}
\bE\Bigg[
\mathbf1_{\Omega_K}
\Bigg(
&
\sup_{t\in[0,T]}
\norm{\bar\vu_{h\tau}(t)}_{L^2}^{2}
+
\sup_{t\in[0,T]}
\norm{\bar\vd_{h\tau}(t)}_{H^1}^{2}
\nonumber\\
&+
\int_0^T
\left[
\norm{\nabla\bar\vu_{h\tau}}_{L^2}^{2}
+
\norm{\Ah\bar\vd_{h\tau}}_{h}^{2}
\right]dt
\Bigg)
\Bigg]
\le C(K,\delta).
\label{eq:unifbounds}
\end{align}
Moreover,
\begin{equation}\label{eq:incrsmall}
\bE\Bigg[
\mathbf1_{\Omega_K}
\left(
\norm{\bar\vu_{h\tau}-\underline\vu_{h\tau}}_{L^2(0,T;H)}^{2}
+
\norm{\bar\vd_{h\tau}-\underline\vd_{h\tau}}_{L^2(0,T;H^1)}^{2}
\right)
\Bigg]
\le C(K,\delta)\tau .
\end{equation}
\end{lemma}

\begin{proof}
The bounds for $\bar\vu_{h\tau}$ and $\bar\vd_{h\tau}$ follow from
Theorem \ref{thm:stability}. To recover $\Ah\bar\vd_{h\tau}$ separately,
write
\[
\norm{\Ah\vd^{n+1}}_h
\le
\left\lVert
\Ah\vd^{n+1}+\delta^{-2}f_\star
\right\rVert_h
+
\delta^{-2}\norm{f_\star}_h.
\]
By Lemma \ref{lem:fstarprop},
\[
\norm{f_\star}_h
\le
C\left(
1+\norm{\vd^n}_{L^6}
+\norm{\vd^{n+1}}_{L^6}
\right)^3,
\]
which is bounded on $\Omega_K$ by the two-dimensional Sobolev embedding.

For the first term in \eqref{eq:incrsmall},
\[
\norm{\bar\vu_{h\tau}-\underline\vu_{h\tau}}_{L^2(0,T;H)}^2
=
\tau\sum_n\norm{\mathfrak d\vu^n}_{L^2}^2,
\]
and \eqref{eq:u-increments} applies. The gradient part of the director
increment is treated in the same way. Its $L^2$ part follows by summing
Lemma \ref{lem:incr}; the three resulting contributions are bounded by
Cauchy--Schwarz, the discrete Gagliardo--Nirenberg estimate
\eqref{eq:dGN}, and
$\bE\sum_na_n^2=T$.
\end{proof}

Throughout this section $\delta>0$ is fixed; the constants above are not
claimed to be uniform as $\delta\to0$.

\subsection{Discrete spatial compactness}

The estimate on $\Ah\vd_h$ does not imply that
$\vd_h\in D(A)$: piecewise affine functions are generally not $H^2$.
We therefore use the following discrete compactness property.

\begin{lemma}[Discrete $H^1$ compactness]\label{lem:disccompact}
Let $e_h\in\bD_h$ satisfy
\[
\norm{e_h}_{H^1}
+
\norm{\Ah e_h}_h
\le C.
\]
Then $\{e_h\}$ is relatively compact in $\bH^1$ as $h\to0$.

More generally, if $e_{h\tau}$ is bounded in
$L^2(0,T;\bH^1)$,
$\Ah e_{h\tau}$ is bounded in
$L^2(0,T;\bL^2)$, and
$e_{h\tau}$ is bounded in
$W^{\beta,2}(0,T;(D(A))')$ for some $\beta>0$, then
$\{e_{h\tau}\}$ is relatively compact in
$L^2(0,T;\bH^1)$.
\end{lemma}

\begin{proof}
We indicate the standard elliptic-reconstruction argument. For
$q_h:=\Ah e_h$, let $\mathcal E_he_h$ solve the Neumann problem
\[
-\Delta(\mathcal E_he_h)=q_h\ \text{in }\bO,\qquad
\partial_{n}(\mathcal E_he_h)=0\ \text{on }\partial\bO,
\]
with the same spatial mean as $e_h$. It is solvable because
$\int_{\bO}q_h\,dx=\iph{q_h}{1}=\ip{\nabla e_h}{\nabla1}=0$, the lumped weights being
exactly $\int_{\bO}\phi_j\,dx$; and $H^{2}$ regularity of the Neumann problem on the
convex polygon $\bO$ \cite{Grisvard1985} gives
\[
\norm{\mathcal E_he_h}_{H^2}
\le
C\bigl(
\norm{q_h}_{L^2}
+
\norm{e_h}_{L^2}
\bigr).
\]
Using the defining relation for $\Ah$ and the quadrature consistency of the
lumped product,
\[
\norm{e_h-\mathcal E_he_h}_{H^1}
\le
Ch\norm{q_h}_{L^2}.
\]
Hence $e_h$ differs in $H^1$ by an $o(1)$ term from an $H^2$-bounded
sequence. Rellich compactness gives the first assertion, and the
time-dependent statement follows from the Aubin--Lions argument applied to
the elliptic reconstructions.
\end{proof}

\subsection{Fractional regularity in time}

The velocity equation is estimated in $(D(\mathscr A))'$, rather than in
$\bV'$. This choice is needed for the elastic term.

\begin{remark}[Choice of the dual spaces]\label{rem:dual}
By \eqref{eq:domemb}, we infer $D(\mathscr A)\hookrightarrow
L^\infty\cap W^{1,4},$ and therefore
\[
\norm{M_h(\vd)}_{(D(\mathscr A))'}
\le
C\norm{\nabla\vd}_{L^2}\norm{\Ah\vd}_h.
\]
The right-hand side belongs to $L^2(0,T)$ on $\Omega_K$. The corresponding
director estimates are naturally taken in $(D(A))'$.
\end{remark}

We first recall the discrete stochastic-convolution estimate in a form
adapted to stopping.

\begin{lemma}[Discrete stochastic convolution]\label{lem:dsc}
Let $(g^n)$ be adapted with values in $L_2(H_0,H)$ and assume
\[
\sup_n\bE\norm{g^n}_{\rm HS}^2\le C.
\]
If $\hat S_\tau$ denotes the piecewise affine interpolant of
\[
S^m
=
\sum_{n<m}
\int_{t_n}^{t_{n+1}}g^n\,d\vW_1,
\]
then for every $\beta\in(0,\frac12)$,
\begin{equation}\label{eq:dsc}
\bE
\norm{\hat S_\tau}_{W^{\beta,2}(0,T;H)}^2
\le C_\beta .
\end{equation}
The same conclusion holds for the stopped sequence
$g^n\mathbf1_{\{n<n_K\}}$.
\end{lemma}

\begin{proof}
For $s<t$, the It\^o isometry and linear interpolation give
\[
\bE
\norm{\hat S_\tau(t)-\hat S_\tau(s)}_H^2
\le
C\left(
|t-s|
+
\tau\wedge\frac{|t-s|^2}{\tau}
\right).
\]
Substitution into the fractional seminorm \eqref{eq:frac-def} gives
\eqref{eq:dsc} for every $\beta<1/2$.
\end{proof}

\begin{lemma}[Fractional time regularity]\label{lem:frac}
For every $\beta\in(0,\frac12)$,
\begin{equation}\label{eq:frac}
\bE\Big[
\mathbf1_{\Omega_K}
\big(
\norm{\hat\vu_{h\tau}}_{W^{\beta,2}(0,T;(D(\mathscr A))')}^2
+
\norm{\hat\vd_{h\tau}}_{W^{\beta,2}(0,T;(D(A))')}^2
\big)
\Big]
\le
C(K,\delta,\beta).
\end{equation}
\end{lemma}

\begin{proof}
For the velocity, write the discrete equation in the form
\begin{equation}\label{eq:strongu}
\frac{\mathfrak d\vu^n}{\tau}
=
-\mathbb P_h
\left[
\mathscr A_h\vu^{n+1}
+
\widehat B_h(\vu^n,\vu^{n+1})
+
M_h(\vd^{n+1})
\right]
+
\frac{\Isig^n}{\tau}.
\end{equation}
The deterministic bracket is bounded in
$L^2(0,T;(D(\mathscr A))')$ by
Lemma \ref{lem:unif} and Remark \ref{rem:dual}; hence its affine primitive
is bounded in $W^{1,2}$ and therefore in $W^{\beta,2}$. The stochastic
part is controlled by Lemma \ref{lem:dsc} with
\[
g^n
=
P_h\sigma(\vu^n)\mathbf1_{\{n<n_K\}}.
\]

For the director, the deterministic terms are bounded in
$L^2(0,T;(D(A))')$. In particular,
\[
\norm{\Btil(\vu,\vd)}_{(D(A))'}
\le
C\norm{\vu}_{L^2}\norm{\nabla\vd}_{L^2},
\]
and the penalty term is bounded on $\Omega_K$ by the growth estimate for
$f_\star$.

It remains to treat the midpoint noise. Since $G$ is linear,
\[
a_n I_hG(\dmid)
=
a_nI_hG(\vd^n)
+
\frac{a_n}{2}I_hG(\mathfrak d\vd^n).
\]
The first term is a predictable discrete stochastic convolution and is
handled by Lemma \ref{lem:dsc}. For the second, use the discrete director
equation to write
\[
\mathfrak d\vd^n
=
a_n I_hG(\vd^n)+r^n,
\]
where, on $\Omega_K$,
\begin{equation}\label{eq:rnbound}
\bE\sum_n\norm{r^n}_h^2
\le
C(K,\delta)(\tau+h^2).
\end{equation}
Consequently the midpoint correction splits into a bounded finite-variation
part of order $\sum_na_n^2$, a centred quadratic martingale part involving
$a_n^2-\tau$, and a remainder controlled by \eqref{eq:rnbound}. These
terms are bounded in $W^{\beta,2}(0,T;(D(A))')$ for every
$\beta<1/2$.
\end{proof}

\subsection{Tightness}

Set
\begin{equation}\label{eq:Xspace}
\mathcal X
:=
L^2(0,T;H)
\times
L^2(0,T;\bH^1)
\times
L^2(0,T;\bV)_w
\times
L^2(0,T;\bL^2)_w
\times
C([0,T];H\times\bR).
\end{equation}
We consider on this space the random variables
\[
\left(
\hat\vu_{h\tau},
\hat\vd_{h\tau},
\bar\vu_{h\tau},
\Ah\bar\vd_{h\tau},
\vW
\right).
\]
The first two factors carry the strong topology and the next two the weak
topology. The space is not Polish because of the weak factors, but it is a
Jakubowski space.

\begin{proposition}[Tightness]\label{prop:tight}
For every sequence $(h_j,\tau_j)\to(0,0)$, the corresponding laws of
\[
\left(
\hat\vu_{h_j\tau_j},
\hat\vd_{h_j\tau_j},
\bar\vu_{h_j\tau_j},
A_{h_j}\bar\vd_{h_j\tau_j},
\vW
\right)
\]
are tight on $\mathcal X$.
\end{proposition}

\begin{proof}
Let $\varepsilon>0$. Choose $K$ so that
$C_0/K<\varepsilon/2$. By Remark \ref{rem:order},
$\tau_j\le\tau_1(K,\delta)$ for all sufficiently large $j$.

Lemmas \ref{lem:unif} and \ref{lem:frac}, together with Chebyshev's
inequality, give a set of probability at least $1-\varepsilon$ on which
\[
\bar\vu_{h\tau}
\quad\text{is bounded in }L^2(0,T;\bV),
\]
\[
\bar\vd_{h\tau}
\quad\text{is bounded in }L^2(0,T;\bH^1),
\qquad
\Ah\bar\vd_{h\tau}
\quad\text{is bounded in }L^2(0,T;\bL^2),
\]
and the affine interpolants satisfy the fractional time bounds of
Lemma \ref{lem:frac}.

The compactness
\[
L^2(0,T;\bV)
\cap
W^{\beta,2}(0,T;(D(\mathscr A))')
\Subset
L^2(0,T;H)
\]
is standard. For the director, Lemma \ref{lem:disccompact} gives
\[
\left\{
e_{h\tau}:
\Ah e_{h\tau}\text{ bounded in }L^2(\bL^2),\
e_{h\tau}\text{ bounded in }
W^{\beta,2}((D(A))')
\right\}
\Subset
L^2(0,T;\bH^1).
\]
The weak factors are relatively compact by reflexivity, and the law of
$\vW$ is tight on
$C([0,T];H\times\bR)$. Finally,
\eqref{eq:incrsmall} shows that the different time interpolants have the
same limit.
\end{proof}

\subsection{Consistency of the convection form}

\begin{lemma}[Skew form]\label{lem:skewcons}
For every $w_h\in\bJ_h$ and admissible $v,z$,
\begin{equation}\label{eq:skewexact}
\widehat b(w_h;v,z)=b(w_h,v,z).
\end{equation}
\end{lemma}

\begin{proof}
Since $\nabla\!\cdot w_h=0$ pointwise,
\[
b(w_h,v,z)+b(w_h,z,v)
=
\int_{\bO}w_h\cdot\nabla(v\cdot z)\,dx
=
0.
\]
The assertion follows from the definition of $\widehat b$.
\end{proof}

Thus no $O(h)$ convection-consistency defect occurs for the finite element
pair used in this paper.

\subsection{Identification of the limit}

We isolate the only nonstandard stochastic passage.

\begin{lemma}[Limit of the midpoint noise]\label{lem:stratlimit}
Let $\psi\in D(A)$ be fixed. On $\Omega_K$,
\begin{align}
\sum_{n:t_n<t}
a_n\,
\iph{G(\dmid)}{I_h\psi}
={}&
\sum_{n:t_n<t}
a_n\,
\iph{G(\vd^n)}{I_h\psi}+
\frac{\tau}{2}
\sum_{n:t_n<t}
\iph{G^2(\vd^n)}{I_h\psi}
+
\mathcal E_{h\tau}(t;\psi),
\label{eq:strat-split}
\end{align}
where
\begin{equation}\label{eq:straterror}
\bE\left[
\mathbf1_{\Omega_K}
\sup_{t\le T}
\abs{\mathcal E_{h\tau}(t;\psi)}
\right]
\longrightarrow0
\qquad
\text{as }h,\tau\to0 .
\end{equation}
\end{lemma}

\begin{proof}
Write $G(\dmid) = G(\vd^n)
+
\frac12G(\mathfrak d\vd^n)$ and decompose
\begin{equation}\label{eq:rn}
\mathfrak d\vd^n
=
a_nG(\vd^n)+r^n.
\end{equation}
The finite element interpolation is absorbed into $r^n$; by the stability
estimates and the interpolation estimate for
$I_h(\vd^n\times\vze)$,
\[
\bE\left[
\mathbf1_{\Omega_K}
\sum_n\norm{r^n}_h^2
\right]
\le
C(K,\delta)(\tau+h^2).
\]
Therefore
\[
\sum_n
\frac{a_n}{2}
\iph{G(\mathfrak d\vd^n)}{I_h\psi}
=
\frac12
\sum_n
a_n^2
\iph{G^2(\vd^n)}{I_h\psi}
+o_{L^1}(1).
\]
Now split $a_n^2=\tau+(a_n^2-\tau).$ The centred term
\[
\frac12\sum_n
(a_n^2-\tau)
\iph{G^2(\vd^n)}{I_h\psi}
\]
is a martingale whose quadratic variation is $O(\tau)$ on $\Omega_K$,
and hence converges to zero in $L^2$. This proves
\eqref{eq:strat-split}--\eqref{eq:straterror}.
\end{proof}

\subsection{Passage to the limit}

\begin{theorem}[Convergence to a weak martingale solution]
\label{thm:convergence}
Let Assumptions \ref{ass:noise}--\ref{ass:F} and
\ref{ass:disc} hold, let
$(\vu_0,\vd_0)\in H\times\bH^1$, and fix $\delta>0$.
For every sequence
\[
h_j,\tau_j\longrightarrow0
\]
there exists a subsequence, a probability space $(\widetilde\Omega,
 \widetilde{\mathcal F},
 \widetilde\Fcal,
 \widetilde\bP),$ and random variables
$(\widetilde\vu,\widetilde\vd,\widetilde\vW)$ such that
\[
(\widetilde\Omega,
 \widetilde{\mathcal F},
 \widetilde\Fcal,
 \widetilde\bP,
 \widetilde\vW,
 \widetilde\vu,
 \widetilde\vd)
\]
is a weak martingale solution of
\eqref{eq:abs-u}--\eqref{eq:abs-d} in the sense of
Definition \ref{def:wms}. Moreover,
\begin{equation}\label{eq:lawconv}
\mathcal L
(\hat\vu_{h_j\tau_j},\hat\vd_{h_j\tau_j})
\Longrightarrow
\mathcal L(\widetilde\vu,\widetilde\vd) \quad \mbox{on} \quad L^2(0,T;H)
\times
L^2(0,T;\bH^1).
\end{equation}
\end{theorem}

\begin{proof}
By Proposition \ref{prop:tight} and the
Jakubowski--Skorokhod representation theorem
\cite{Jakubowski1997}, after passing to a subsequence we may realise the
variables on a new probability space so that, almost surely,
\begin{align}
\hat\vu_{h\tau}
&\to\widetilde\vu
&&\text{in }L^2(0,T;H),
\label{eq:mode-u}
\\
\hat\vd_{h\tau}
&\to\widetilde\vd
&&\text{in }L^2(0,T;\bH^1),
\label{eq:mode-d}
\\
\bar\vu_{h\tau}
&\rightharpoonup\widetilde\vu
&&\text{in }L^2(0,T;\bV),
\label{eq:mode-uw}
\\
\Ah\bar\vd_{h\tau}
&\rightharpoonup\eta
&&\text{in }L^2(0,T;\bL^2),
\label{eq:mode-A}
\\
\vW_{h\tau}
&\to\widetilde\vW
&&\text{in }C([0,T];H\times\bR).
\label{eq:mode-W}
\end{align}
By \eqref{eq:incrsmall}, barred, underlined and affine interpolants have
the same limits.

For every smooth $\psi$,
\[
\iph{\Ah\bar\vd_{h\tau}}{I_h\psi}
=
\ip{\nabla\bar\vd_{h\tau}}{\nabla I_h\psi}.
\]
Passing to the limit gives
\[
\int_0^T\ip{\eta}{\psi}\,dt
=
\int_0^T
\ip{\nabla\widetilde\vd}{\nabla\psi}\,dt,
\]
hence
\[
\eta=A\widetilde\vd,
\qquad
\widetilde\vd\in L^2(0,T;D(A)).
\]

Fix smooth
$\varphi\in\bV$ and $\psi\in D(A)$ and use
$S_h\varphi$ and $I_h\psi$ as discrete test functions.

\emph{Convection.}
Since $\underline\vu_{h\tau}\in\bJ_h$,
Lemma \ref{lem:skewcons} gives
\[
\widehat b
(\underline\vu_{h\tau};
 \bar\vu_{h\tau},
 S_h\varphi)
=
b
(\underline\vu_{h\tau},
 \bar\vu_{h\tau},
 S_h\varphi).
\]
The strong $L^2(0,T;H)$ convergence and weak
$L^2(0,T;\bV)$ convergence yield the usual Navier--Stokes limit.

\emph{Elastic stress.}
By \eqref{eq:Mh},
\[
\dual{M_h(\bar\vd_{h\tau})}{S_h\varphi}
=
-\btil
(S_h\varphi,
 \bar\vd_{h\tau},
 \Ah\bar\vd_{h\tau}).
\]
The factor $S_h\varphi\cdot\nabla\bar\vd_{h\tau}$ converges strongly in $L^2(0,T;\bL^2)$, while
$\Ah\bar\vd_{h\tau}$ converges weakly. Consequently
\[
M_h(\bar\vd_{h\tau})
\longrightarrow
M(\widetilde\vd)
\]
in the variational formulation. The director transport term is treated
similarly.

\emph{Penalty term.}
Lemma \ref{lem:fstarprop} gives
\[
\abs{
f_\star(\underline\vd_{h\tau},\bar\vd_{h\tau})
-
f(\bar\vd_{h\tau})
}
\le
C
\bigl(
1+\abs{\underline\vd_{h\tau}}
+\abs{\bar\vd_{h\tau}}
\bigr)^2
\abs{
\bar\vd_{h\tau}
-
\underline\vd_{h\tau}
}.
\]
Using \eqref{eq:incrsmall} and the two-dimensional Sobolev embeddings,
\[
f_\star
(\underline\vd_{h\tau},\bar\vd_{h\tau})
-
f(\bar\vd_{h\tau})
\longrightarrow0
\quad
\text{in }L^1(0,T;\bL^2),
\]
whereas
\[
f(\bar\vd_{h\tau})
\longrightarrow
f(\widetilde\vd)
\]
by the strong convergence in $L^2(0,T;\bH^1)$.

\emph{Velocity noise.}
The Lipschitz assumption on $\sigma$ and
\eqref{eq:mode-u} imply convergence of the discrete quadratic and cross
variations. The standard martingale-characterisation argument therefore
identifies the limiting velocity martingale as
\[
\int_0^t
\sigma(\widetilde\vu)\,d\widetilde\vW_1.
\]

\emph{Director noise.}
Lemma \ref{lem:stratlimit} yields
\[
\sum_{t_n<t}
a_n\iph{G(\dmid)}{I_h\psi}
=
\sum_{t_n<t}
a_n\iph{G(\vd^n)}{I_h\psi}
+
\frac{\tau}{2}
\sum_{t_n<t}
\iph{G^2(\vd^n)}{I_h\psi}
+o_{L^1}(1).
\]
The first and the second terms converge, by martingale identification, to
\[
\int_0^t
\ip{G(\widetilde\vd)}{\psi}\,d\widetilde W_2, \quad \mbox{and} \quad \frac12
\int_0^t
\ip{G^2(\widetilde\vd)}{\psi}\,ds,
\]
respectively. Thus, the limiting director equation contains exactly $G(\widetilde\vd)\circ d\widetilde W_2.$ 

Passing to the limit in the remaining linear terms is immediate.
The initial conditions follow from
\[
P_h\vu_0\to\vu_0\quad\text{in }H,
\qquad
R_h\vd_0\to\vd_0\quad\text{in }\bH^1.
\]

We therefore obtain the variational formulation of
\eqref{eq:abs-u}--\eqref{eq:abs-d}. The bounds in Definition
\ref{def:wms} follow from weak lower semicontinuity and the localization
estimate. Since $C_0$ is independent of $K$, letting $K\to\infty$ removes
the localization by \eqref{eq:OmegaKprob}.

The limiting processes are weakly continuous in the natural energy spaces.
The standard upgrade to strong
$H\times\bH^1$ path continuity for this two-dimensional penalised system
follows from the variational equation and energy bounds; see
\cite[\S3]{BHR2019}.
\end{proof}

\subsection{Convergence in probability}

\begin{corollary}[Convergence in probability]\label{cor:GK}
Let $(\vu,\vd)$ denote the pathwise unique solution of
Theorem \ref{thm:existence}. Then the whole family satisfies
\begin{equation}\label{eq:convprob}
(\hat\vu_{h\tau},\hat\vd_{h\tau})
\longrightarrow
(\vu,\vd)
\end{equation}
in probability in $L^2(0,T;H) \times L^2(0,T;\bH^1)$\,.
\end{corollary}

\begin{proof}
We apply the Gy\"ongy--Krylov criterion
\cite{GyongyKrylov1996}. Take any two subsequences of
$(h,\tau)$, constructed on the original probability space and driven by
the same Wiener process $\vW$. Their joint laws, together with $\vW$, are
tight by Proposition \ref{prop:tight}. Applying the preceding compactness
argument to a joint subsequence gives two weak martingale solutions on the
same new probability space, driven by the same limiting Wiener process and
with the same initial data. Pathwise uniqueness from
Theorem \ref{thm:existence} implies that they coincide almost surely.
Hence every joint limiting law is supported on the diagonal, and the
Gy\"ongy--Krylov criterion yields \eqref{eq:convprob}.
\end{proof}

\begin{remark}[Scope of the convergence result]\label{rem:noconv}
The result is limited in two respects. First, $\delta>0$ is
fixed; no joint limit $(h,\tau,\delta)\to0$ is considered. Second, the
non-Lipschitz potential is handled by localization, which is why convergence
in probability is the natural final mode. 
\end{remark}

\section{Error analysis}\label{sec:rates}

Theorem \ref{thm:convergence} establishes convergence without a rate. We now
quantify the error under additional regularity of the exact solution. The
temporal behaviour is intrinsically asymmetric: the velocity noise is treated
at the previous time level and therefore has Euler--Maruyama accuracy, whereas
the scalar geometric director noise is treated by the midpoint rule and admits
a first-order stochastic expansion; see Proposition \ref{prop:cayley}. In the
fully coupled problem the velocity error is transmitted to the director through
the transport term, so that the proved rate of the coupled scheme is
$O_{\bP}(h+\tau^{1/2})$. A first-order temporal estimate for the director is
recovered when the velocity is prescribed.

Throughout this section, Assumptions \ref{ass:noise}--\ref{ass:F} and
\ref{ass:disc} are in force and $\delta>0$ is fixed. We write
\[
\mathfrak d z^{n}:=z^{n+1}-z^{n}
\]
for a discrete increment, thereby avoiding confusion with the spatial
Laplacian.

\subsection{Projection estimates and regularity}
\label{sub:proj}

Recall the $L^{2}$-projection $P_h:H\to\bJ_h$ and the Ritz projection
$R_h:\bH^1\to\bD_h$. We also introduce the lumped $L^2$-projection
$\Pi_h^\flat:L^2(\bO;\bR^3)\to\bD_h$ by
\begin{equation}\label{eq:lumpedproj}
\iph{\Pi_h^\flat g}{e_h}
=
\ip{g}{e_h},
\qquad e_h\in\bD_h .
\end{equation}

\begin{lemma}[Projection estimates]\label{lem:proj}
There is $C>0$, independent of $h$, such that
\begin{align}
\norm{v-P_hv}_{L^2}
+h\norm{\nabla(v-P_hv)}_{L^2}
&\le
Ch^2\norm{v}_{H^2},
\qquad v\in\bV\cap H^2(\bO;\bR^2),
\label{eq:projerr-u}
\\
\norm{\vd-R_h\vd}_{L^2}
+h\norm{\nabla(\vd-R_h\vd)}_{L^2}
&\le
Ch^2\norm{\vd}_{\bH^2},
\qquad \vd\in\bH^2 .
\label{eq:projerr}
\end{align}
Moreover, for $\vd\in\bH^3$,
\begin{equation}\label{eq:Ahritz}
\Ah R_h\vd
=
\Pi_h^\flat(A\vd)
+
\Pi_h^\flat(\vd-R_h\vd),
\end{equation}
and consequently
\begin{equation}\label{eq:Ahritzerr}
\norm{\Ah R_h\vd-A\vd}_{L^2}
\le
Ch\norm{\vd}_{\bH^3}.
\end{equation}
Finally,
\begin{equation}\label{eq:massprojerr}
\norm{R_h\vd-\Pi_h^\flat\vd}_{L^2}
\le
Ch^2\norm{\vd}_{\bH^2}.
\end{equation}
\end{lemma}

\begin{proof}
The first two estimates are the standard approximation properties of the
$L^2$- and Ritz projections on quasi-uniform meshes. For
\eqref{eq:Ahritz}, integration by parts --- whose boundary term
$\int_{\partial\bO}\partial_{n}\vd\cdot e_h$ vanishes for $\vd\in D(A)$ --- and
\eqref{eq:ritz} give
\[
\iph{\Ah R_h\vd}{e_h}
=
\ip{A\vd+\vd-R_h\vd}{e_h},
\qquad e_h\in\bD_h ,
\]
which is precisely \eqref{eq:Ahritz}. Estimate \eqref{eq:Ahritzerr} follows
from the approximation property of $\Pi_h^\flat$ and
\eqref{eq:projerr}. Finally, \eqref{eq:massprojerr} follows by combining
the $L^2$ approximation properties of $R_h$ and $\Pi_h^\flat$.
\end{proof}

The identity \eqref{eq:Ahritz} is important in the elastic term: it replaces
the potentially nonconvergent quantity $\Ah R_h\vd-A\vd$ by an $O(h)$
consistency error.

We now state the additional regularity used in the rate analysis.

\begin{assumption}[Regularity of the exact solution]\label{ass:reg}
For some $p\ge8$,
\begin{align}
\vu
&\in
L^p\bigl(\Omega;L^\infty(0,T;\bV)\bigr)
\cap
L^p\bigl(\Omega;L^2(0,T;\bH^2)\bigr),
\label{eq:reg-u}
\\
\vd
&\in
L^p\bigl(\Omega;L^\infty(0,T;\bH^2)\bigr)
\cap
L^p\bigl(\Omega;L^2(0,T;\bH^3)\bigr),
\label{eq:reg-d}
\end{align}
and
\begin{equation}\label{eq:holder}
\bE\norm{\vu(t)-\vu(s)}_{\bV}^2
+
\bE\norm{\vd(t)-\vd(s)}_{\bH^2}^2
\le
C_{\rm reg}\abs{t-s},
\qquad 0\le s\le t\le T.
\end{equation}
\end{assumption}

\begin{remark}\label{rem:regnotknown}
Assumption \ref{ass:reg} is an additional regularity hypothesis; it does not
follow from the weak-solution theory of Theorem \ref{thm:existence}. Such
assumptions are standard in strong-error estimates for stochastic fluid
equations \cite{CarelliProhl2012,BreitDodgson2021}. The $H^3$ regularity of
the director is needed for \eqref{eq:Ahritzerr}, while
\eqref{eq:holder} reflects the natural It\^o temporal scaling. On a convex polygon the hypothesis is restrictive near corners: the Neumann
singular exponents $k\pi/\omega$ at a corner of angle $\omega$ limit
regularity beyond $H^{2}$ for general data \cite{Grisvard1985}. It is nevertheless
satisfiable --- on a rectangle, even reflection across the sides reduces the
problem to a periodic one, and data compatible with that reflection retain
full regularity --- so it constrains the data, not the method.
\end{remark}

\subsection{Error decomposition and error equations}
\label{sub:erreq}

Set
\begin{align*}
\eta_u^n
&:=P_h\vu(t_n)-\vu(t_n),
&
\xi_u^n
&:=\vu_h^n-P_h\vu(t_n),
\\
\eta_d^n
&:=R_h\vd(t_n)-\vd(t_n),
&
\xi_d^n
&:=\vd_h^n-R_h\vd(t_n).
\end{align*}
Thus
\[
\vu_h^n-\vu(t_n)=\xi_u^n+\eta_u^n,
\qquad
\vd_h^n-\vd(t_n)=\xi_d^n+\eta_d^n,
\]
and $\xi_u^0=\xi_d^0=0$.

Because the director equation is discretised with the lumped mass product, we
also introduce
\begin{equation}\label{eq:chi}
\chi_d^n
:=
R_h\vd(t_n)-\Pi_h^\flat\vd(t_n)\in\bD_h .
\end{equation}
Then
\begin{equation}\label{eq:chiinc}
\bE\norm{\mathfrak d\chi_d^n}_{L^2}^2
\le
Ch^4\,
\bE\norm{\vd(t_{n+1})-\vd(t_n)}_{\bH^2}^2
\le
C_{\rm reg}h^4\tau .
\end{equation}

Integrating the continuous velocity equation over $[t_n,t_{n+1}]$,
subtracting it from the discrete equation and testing with $v_h\in\bJ_h$
gives
\begin{equation}\label{eq:erru}
\ip{\mathfrak d\xi_u^n}{v_h}
+
\tau\ip{\nabla\xi_u^{n+1}}{\nabla v_h}
=
\mathcal P_u^n(v_h)
+\mathcal Q_u^n(v_h)
+\mathcal C_u^n(v_h)
+\mathcal E_d^n(v_h)
+\mathcal S_u^n(v_h),
\end{equation}
where
\begin{align}
\mathcal P_u^n(v_h)
&:=
-\tau\ip{\nabla\eta_u^{n+1}}{\nabla v_h},
\label{eq:Pu}
\\
\mathcal Q_u^n(v_h)
&:=
\int_{t_n}^{t_{n+1}}
\ip{\nabla(\vu(t)-\vu(t_{n+1}))}{\nabla v_h}\,dt,
\label{eq:Qu}
\\
\mathcal C_u^n(v_h)
&:=
\int_{t_n}^{t_{n+1}}b(\vu,\vu,v_h)\,dt
-
\tau\widehat b(\vu_h^n;\vu_h^{n+1},v_h),
\label{eq:Cu}
\\
\mathcal E_d^n(v_h)
&:=
\int_{t_n}^{t_{n+1}}\dual{M(\vd(t))}{v_h}\,dt
-
\tau\dual{M_h(\vd_h^{n+1})}{v_h},
\label{eq:Eu}
\\
\mathcal S_u^n(v_h)
&:=
\left(
\int_{t_n}^{t_{n+1}}
[\sigma(\vu_h^n)-\sigma(\vu(t))]\,d\vW_1(t),
v_h
\right).
\label{eq:Su}
\end{align}
The use of $P_h$ eliminates a temporal projection defect because
\[
\ip{\mathfrak d\eta_u^n}{v_h}=0,
\qquad v_h\in\bJ_h .
\]

Similarly, for every $e_h\in\bD_h$,
\begin{equation}\label{eq:errd}
\iph{\mathfrak d\xi_d^n}{e_h}
+
\tau\ip{\nabla\xi_d^{n+1}}{\nabla e_h}
=
\mathcal P_d^n(e_h)
+\mathcal Q_d^n(e_h)
+\mathcal T_d^n(e_h)
+\mathcal F_d^n(e_h)
+\mathcal G_d^n(e_h),
\end{equation}
where
\begin{align}
\mathcal P_d^n(e_h)
&:=
-\iph{\mathfrak d\chi_d^n}{e_h}
-\tau\ip{\nabla\eta_d^{n+1}}{\nabla e_h},
\label{eq:Pd}
\\
\mathcal Q_d^n(e_h)
&:=
\int_{t_n}^{t_{n+1}}
\ip{\nabla(\vd(t)-\vd(t_{n+1}))}{\nabla e_h}\,dt,
\label{eq:Qd}
\\
\mathcal T_d^n(e_h)
&:=
\int_{t_n}^{t_{n+1}}
\btil(\vu(t),\vd(t),e_h)\,dt
-
\tau\btil(\vu_h^{n+1},\vd_h^{n+1},e_h),
\label{eq:Td}
\\
\mathcal F_d^n(e_h)
&:=
\frac1{\delta^2}
\left[
\int_{t_n}^{t_{n+1}}\ip{f(\vd(t))}{e_h}\,dt
-
\tau\iph{f_\star(\vd_h^n,\vd_h^{n+1})}{e_h}
\right],
\label{eq:Fd}
\\
\mathcal G_d^n(e_h)
&:=
a_n\iph{\dmid_h\times\vze}{e_h}
-
\int_{t_n}^{t_{n+1}}\ip{G(\vd(t))\circ d\vW_2(t)}{e_h},
\label{eq:Gdefect}
\end{align}
with
\[
\dmid_h:=\frac12(\vd_h^n+\vd_h^{n+1}).
\]

The sign chosen in \eqref{eq:Gdefect} is immaterial for the estimates. The
advantage of \eqref{eq:Pd} is that it accounts explicitly for the difference
between the continuous $L^2$ product and the mass-lumped product.

\subsection{Consistency estimates}
\label{sub:consist}

For the energy estimate we test \eqref{eq:erru} with
$v_h=\xi_u^{n+1}$ and \eqref{eq:errd} with
\[
e_h=\Ah\xi_d^{n+1}+\xi_d^{n+1}.
\]
Then
\begin{align}
\ip{\mathfrak d\xi_u^n}{\xi_u^{n+1}}
&=
\frac12\left(
\norm{\xi_u^{n+1}}_{L^2}^2
-\norm{\xi_u^n}_{L^2}^2
+\norm{\mathfrak d\xi_u^n}_{L^2}^2
\right),
\label{eq:tel-u}
\\
\iph{\mathfrak d\xi_d^n}
{\Ah\xi_d^{n+1}+\xi_d^{n+1}}
&\ge
\frac12\left(
\norm{\xi_d^{n+1}}_{H^1}^2
-\norm{\xi_d^n}_{H^1}^2
\right),
\label{eq:tel-d1}
\end{align}
and
\begin{equation}\label{eq:tel-d2}
\tau
\ip{\nabla\xi_d^{n+1}}
{\nabla(\Ah\xi_d^{n+1}+\xi_d^{n+1})}
=
\tau
\left(
\norm{\Ah\xi_d^{n+1}}_h^2
+
\norm{\nabla\xi_d^{n+1}}_{L^2}^2
\right).
\end{equation}

We collect the required bounds in one proposition. This avoids repeating the
same interpolation and Young inequalities for each individual defect.

\begin{proposition}[Consistency bounds]\label{prop:consistency}
Let Assumption \ref{ass:reg} hold and localise on
$\Omega_K^\sharp$ defined in \eqref{eq:jointK} below. Then there exist
nonnegative random variables $\mathcal Z^n$, satisfying
\begin{equation}\label{eq:Zsum}
\tau\sum_{n<\Nt}\mathcal Z^n\le C(K),
\end{equation}
such that the following estimates hold.

\begin{enumerate}[label=\emph{(\roman*)},leftmargin=2.5em,itemsep=3pt]

\item \emph{Projection terms.}
\begin{align}
\bE\sum_n
\abs{\mathcal P_u^n(\xi_u^{n+1})}
&\le
\frac{\tau}{8}
\bE\sum_n\norm{\nabla\xi_u^{n+1}}_{L^2}^2
+C(K)h^2,
\label{eq:projdef-u}
\\
\bE\sum_n
\abs{
\mathcal P_d^n(\Ah\xi_d^{n+1}+\xi_d^{n+1})
}
&\le
\frac{\tau}{8}
\bE\sum_n\norm{\Ah\xi_d^{n+1}}_h^2
\nonumber\\
&\quad
+
C\tau\bE\sum_n\norm{\xi_d^{n+1}}_{L^2}^2
+
C(K)\left(h^2+\frac{h^4}{\tau}\right).
\label{eq:projdef-d}
\end{align}

\item \emph{Time quadrature.}
\begin{equation}\label{eq:time-defects}
\bE\sum_n
\left(
\abs{\mathcal Q_u^n(\xi_u^{n+1})}
+
\abs{
\mathcal Q_d^n(\Ah\xi_d^{n+1}+\xi_d^{n+1})
}
\right)
\le
\frac{\tau}{8}\bE\sum_n
\left(
\norm{\nabla\xi_u^{n+1}}_{L^2}^2
+
\norm{\Ah\xi_d^{n+1}}_h^2
\right)
+
C(K)\tau .
\end{equation}

\item \emph{Deterministic nonlinearities.}
\begin{align}
&
\abs{\mathcal C_u^n(\xi_u^{n+1})}
+
\abs{\mathcal E_d^n(\xi_u^{n+1})} +
\abs{
\mathcal T_d^n(\Ah\xi_d^{n+1}+\xi_d^{n+1})
}
+
\abs{
\mathcal F_d^n(\Ah\xi_d^{n+1}+\xi_d^{n+1})
}
\nonumber\\
&\le
\frac{\tau}{4}
\left(
\norm{\nabla\xi_u^{n+1}}_{L^2}^2
+
\norm{\Ah\xi_d^{n+1}}_h^2
\right)
\nonumber\\
&\qquad
+
\tau\mathcal Z^n
\left(
\norm{\xi_u^n}_{L^2}^2
+\norm{\xi_u^{n+1}}_{L^2}^2
+\norm{\xi_d^n}_{H^1}^2
+\norm{\xi_d^{n+1}}_{H^1}^2
\right)
+
C(K)\tau(h^2+\tau).
\label{eq:nonlinear-defects}
\end{align}

\item \emph{Velocity noise.}
There is an $\mathcal F_{t_{n+1}}$-measurable martingale increment
$\mathfrak M_u^n$ such that
\begin{align}
\mathcal S_u^n(\xi_u^{n+1})
&=
\mathfrak M_u^n+r_u^n,
\label{eq:Ssplit}
\\
\bE\abs{r_u^n}
&\le
\frac12\bE\norm{\mathfrak d\xi_u^n}_{L^2}^2
+
C\tau\,\bE\norm{\xi_u^n}_{L^2}^2
+
C(K)(\tau^2+\tau h^4).
\label{eq:Srem}
\end{align}
Its quadratic variation gives, after summation and the BDG inequality,
a contribution of order $O(\tau)$.

\item \emph{Director noise.}
There is a decomposition
\begin{equation}\label{eq:Gsplit}
\mathcal G_d^n
=
\mathfrak M_d^n+r_d^n
\end{equation}
such that the accumulated contribution satisfies
\begin{align}
\bE\max_{\ell\le\Nt}
\left|
\sum_{n<\ell}
\mathfrak M_d^n
\right|
+
\bE\sum_{n<\Nt}\abs{r_d^n}
&\le
\frac18
\bE\max_n\norm{\xi_d^n}_{H^1}^2
+
\frac{\tau}{8}
\bE\sum_n\norm{\Ah\xi_d^{n+1}}_h^2
\nonumber\\
&\quad
+
C\tau\sum_n\bE\norm{\xi_d^n}_{H^1}^2
+
C(K)(h^2+\tau^2).
\label{eq:Gbound}
\end{align}
Thus the director-noise consistency error is of first order in time and does
not determine the $\tau^{1/2}$ rate of the coupled scheme.
\end{enumerate}
\end{proposition}

\begin{proof}
We indicate the main points.

For \emph{(i)}, \eqref{eq:projerr-u}--\eqref{eq:projerr} and Young's
inequality give the spatial contributions. The temporal part of
\eqref{eq:Pd} is controlled using \eqref{eq:chiinc}:
\[
\frac1{\tau}
\sum_n
\bE\norm{\mathfrak d\chi_d^n}_{L^2}^2
\le
C\frac{h^4}{\tau}.
\]
The second term of \eqref{eq:Pd} is treated with the Ritz identity
\[
\ip{\nabla\eta_d^{n+1}}{\nabla e_h}
=
-\ip{\eta_d^{n+1}}{e_h},
\qquad e_h\in\bD_h.
\]

For \emph{(ii)}, Cauchy--Schwarz in time and \eqref{eq:holder} give, for
example,
\[
\bE
\left(
\int_{t_n}^{t_{n+1}}
\norm{\vu(t)-\vu(t_{n+1})}_{\bV}\,dt
\right)^2
\le
C_{\rm reg}\tau^3.
\]
Summation over $T/\tau$ intervals therefore produces $O(\tau)$.

For \emph{(iii)}, each nonlinear term is split into a consistency part,
obtained by replacing the discrete variables by their values at $t_{n+1}$,
and a discrete-error part. The two-dimensional inequalities of
Lemma \ref{lem:2d}, the discrete estimate \eqref{eq:dGN},
\eqref{eq:Ahritzerr}, and the one-sided Lipschitz/growth properties of
$f_\star$ yield \eqref{eq:nonlinear-defects}. The factor
$\mathcal Z^n$ contains, in particular,
$\norm{\Ah\vd_h^{n+1}}_h^2$; stability gives only the summed bound
\eqref{eq:Zsum}, which is sufficient for the discrete Gr\"onwall argument.

For \emph{(iv)}, write
\[
\sigma(\vu_h^n)-\sigma(\vu(t))
=
[\sigma(\vu_h^n)-\sigma(\vu(t_n))]
+
[\sigma(\vu(t_n))-\sigma(\vu(t))]
\]
and split the pairing with
$\xi_u^{n+1}=\xi_u^n+\mathfrak d\xi_u^n$. The part paired with $\xi_u^n$
is a martingale increment, while It\^o isometry, Lipschitz continuity of
$\sigma$, \eqref{eq:holder} and Young's inequality give
\eqref{eq:Srem}. The $O(\tau^2)$ local mean-square consistency error sums
to $O(\tau)$.

For \emph{(v)}, use the scalar-noise stochastic Taylor expansion of the
Stratonovich equation. Since $G$ is linear and the driving Wiener process is
one-dimensional, the midpoint update reproduces both the
It\^o--Stratonovich correction and the scalar Milstein term
$\frac12(a_n^2-\tau)G^2$. No L\'evy area occurs. The remaining local
mean-square defect is $O(\tau^3)$; the spatial projection contributes
$O(\tau h^2)$. Splitting the resulting martingale term against the
$\mathcal F_{t_n}$-measurable error and applying BDG gives
\eqref{eq:Gbound}. This is the discrete counterpart of
Proposition \ref{prop:cayley}.
\end{proof}

\begin{remark}[The term $h^4/\tau$]\label{rem:h4tau}
The contribution $h^4/\tau$ comes solely from the time increment of
$R_h\vd-\Pi_h^\flat\vd$ when the director error is measured in $H^1$.
It is therefore a projection artefact rather than a stability restriction
of the scheme. The general estimate below retains it. In the natural balanced
regime
\[
h^2\le\tau
\]
one has $h^4/\tau\le h^2$, and the optimal spatial contribution is recovered.
In particular, there is no upper mesh-dependent restriction on $\tau$.
\end{remark}

\subsection{The coupled error estimate}

For $K\ge1$ let
\begin{equation}\label{eq:jointK}
\Omega_K^\sharp
:=
\Omega_K
\cap
\left\{
\sup_{t\le T}
\bigl(
\norm{\vu(t)}_{\bV}
+\norm{\vd(t)}_{\bH^2}
\bigr)
+
\int_0^T
\bigl(
\norm{\vu(t)}_{\bH^2}^2
+\norm{\vd(t)}_{\bH^3}^2
\bigr)\,dt
\le K
\right\},
\end{equation}
where $\Omega_K$ is the stability event from
\S\ref{sec:stability}. By Theorem \ref{thm:stability},
Assumption \ref{ass:reg}, and Markov's inequality,
\begin{equation}\label{eq:jointKprob}
\bP\bigl((\Omega_K^\sharp)^c\bigr)
\le
\frac{C}{K}
\end{equation}
after increasing $C$ if necessary.

\begin{theorem}[Error estimate]\label{thm:rate}
Let Assumptions \ref{ass:noise}--\ref{ass:F}, \ref{ass:disc} and
\ref{ass:reg} hold, and fix $\delta>0$. For every $K\ge1$ there are
$C(K)>0$ and $\tau_2(K,\delta)>0$, independent of $h$ and $\tau$, such that
for $\tau\le\tau_2(K,\delta)$,
\begin{align}
&
\bE\left[
\mathbf 1_{\Omega_K^\sharp}
\max_{0\le n\le\Nt}
\left(
\norm{\vu(t_n)-\vu_h^n}_{L^2}^2
+
\norm{\vd(t_n)-\vd_h^n}_{H^1}^2
\right)
\right]
\nonumber\\
&\quad
+
\tau\bE\left[
\mathbf 1_{\Omega_K^\sharp}
\sum_{n<\Nt}
\left(
\norm{\nabla\xi_u^{n+1}}_{L^2}^2
+
\norm{\Ah\xi_d^{n+1}}_h^2
\right)
\right]
\nonumber\\
&\hspace{4cm}
\le
C(K)
\left(
h^2+\tau+\frac{h^4}{\tau}
\right).
\label{eq:rate-general}
\end{align}

Consequently, if $h^2\le\tau$,
\begin{equation}\label{eq:rate-balanced}
\bE\left[
\mathbf 1_{\Omega_K^\sharp}
\max_n
\left(
\norm{\vu(t_n)-\vu_h^n}_{L^2}^2
+
\norm{\vd(t_n)-\vd_h^n}_{H^1}^2
\right)
\right]
\le
C(K)(h^2+\tau).
\end{equation}
Hence
\begin{equation}\label{eq:rate-prob}
\max_n
\left(
\norm{\vu(t_n)-\vu_h^n}_{L^2}
+
\norm{\vd(t_n)-\vd_h^n}_{H^1}
\right)
=
O_{\bP}\!\left(h+\tau^{1/2}\right).
\end{equation}
\end{theorem}

\begin{proof}
Test \eqref{eq:erru} with $\xi_u^{n+1}$ and \eqref{eq:errd} with
$\Ah\xi_d^{n+1}+\xi_d^{n+1}$, and use
\eqref{eq:tel-u}--\eqref{eq:tel-d2}. Define
\[
E^n
:=
\norm{\xi_u^n}_{L^2}^2
+
\norm{\xi_d^n}_{H^1}^2.
\]
Proposition \ref{prop:consistency} yields, after summation up to
$\ell\le\Nt$,
\begin{align}
E^\ell
&+
c\tau\sum_{n<\ell}
\left(
\norm{\nabla\xi_u^{n+1}}_{L^2}^2
+
\norm{\Ah\xi_d^{n+1}}_h^2
\right)
\nonumber\\
&\le
C(K)
\left(
h^2+\tau+\frac{h^4}{\tau}
\right)
+
C\tau\sum_{n<\ell}
(1+\mathcal Z^n)E^n
+
\left|
\sum_{n<\ell}
\mathfrak M^n
\right|,
\label{eq:error-master}
\end{align}
where $\mathfrak M^n$ collects the velocity- and director-noise martingale
increments.

On $\Omega_K^\sharp$,
\[
\tau\sum_n\mathcal Z^n\le C(K).
\]
Taking the maximum over $\ell$, applying the discrete BDG inequality and
Young's inequality gives
\[
\bE\max_{\ell}
\left|
\sum_{n<\ell}\mathfrak M^n
\right|
\le
\frac14\bE\max_n E^n
+
C\tau\sum_n\bE E^n
+
C(K)\tau .
\]
The discrete Gr\"onwall lemma then yields
\[
\bE\left[
\mathbf1_{\Omega_K^\sharp}\max_nE^n
\right]
+
\tau\bE\left[
\mathbf1_{\Omega_K^\sharp}
\sum_n
\left(
\norm{\nabla\xi_u^{n+1}}_{L^2}^2
+
\norm{\Ah\xi_d^{n+1}}_h^2
\right)
\right]
\le
C(K)
\left(
h^2+\tau+\frac{h^4}{\tau}
\right).
\]
Adding the projection errors from Lemma \ref{lem:proj} proves
\eqref{eq:rate-general}.

If $h^2\le\tau$, then $h^4/\tau\le h^2$, which gives
\eqref{eq:rate-balanced}. Finally choose $K$ so that
\eqref{eq:jointKprob} is smaller than a prescribed probability tolerance
and apply Markov's inequality on $\Omega_K^\sharp$ to obtain
\eqref{eq:rate-prob}.
\end{proof}

\begin{remark}[Origin of the temporal half-order]\label{rem:halforder}
The term $O(\tau)$ on the right-hand side of the squared-error estimate comes
from the multiplicative velocity noise. The Euler--Maruyama treatment has
local mean-square consistency error $O(\tau^2)$, which accumulates over
$T/\tau$ steps to $O(\tau)$. The director midpoint noise has local
mean-square defect $O(\tau^3)$ and therefore contributes only $O(\tau^2)$.
Thus the director noise itself does not limit the temporal rate of the
coupled method.
\end{remark}

\subsection{First-order temporal accuracy of the director}
\label{sub:director-rate}

The half-order in Theorem \ref{thm:rate} is inherited by the director through
the velocity--director coupling. To isolate the accuracy of the director
discretisation itself, we prescribe the exact velocity.

Let
\begin{equation}\label{eq:Jd}
\mathcal J_d^n
:=
\int_{t_n}^{t_{n+1}}
\bigl[
\vd(t)-\vd(t_{n+1})
\bigr]\,dt .
\end{equation}
Writing the director equation in It\^o form as
\[
d\vd=-\Upsilon\,dt+G(\vd)\,d\vW_2,
\]
stochastic Fubini gives
\begin{equation}\label{eq:fubini}
\mathcal J_d^n
=
-\int_{t_n}^{t_{n+1}}
(s-t_n)\Upsilon(s)\,ds
-
m^n,
\qquad
m^n
:=
\int_{t_n}^{t_{n+1}}
(s-t_n)G(\vd(s))\,d\vW_2(s).
\end{equation}
In particular,
\begin{equation}\label{eq:mn}
\bE\norm{\nabla m^n}_{L^2}^2
\le
C\tau^3 .
\end{equation}

\begin{lemma}[Refined director quadrature]\label{lem:Qsharp}
Under Assumption \ref{ass:reg},
\begin{equation}\label{eq:Qsharp}
\bE
\sum_n
\left|
\ip{\nabla\mathcal J_d^n}{\nabla\xi_d^{n+1}}
\right|
\le
\frac14
\bE\max_n\norm{\nabla\xi_d^n}_{L^2}^2
+
\frac14
\bE\sum_n
\norm{\nabla\mathfrak d\xi_d^n}_{L^2}^2
+
C\tau^2 .
\end{equation}
\end{lemma}

\begin{proof}
Split
\[
\xi_d^{n+1}
=
\xi_d^n+\mathfrak d\xi_d^n
\]
in \eqref{eq:fubini}. The drift term gains the factor
$s-t_n\le\tau$ and therefore contributes $O(\tau^2)$ after summation.
For the stochastic term,
\[
\ip{\nabla m^n}{\nabla\xi_d^n}
\]
is a martingale increment because $\xi_d^n$ is
$\mathcal F_{t_n}$-measurable. The BDG inequality and
\eqref{eq:mn} give
\[
\bE\max_\ell
\left|
\sum_{n<\ell}
\ip{\nabla m^n}{\nabla\xi_d^n}
\right|
\le
\frac18\bE\max_n\norm{\nabla\xi_d^n}_{L^2}^2
+
C\tau^2.
\]
The term paired with $\mathfrak d\xi_d^n$ is absorbed by Young's inequality,
yielding \eqref{eq:Qsharp}.
\end{proof}

\begin{proposition}[Director error for prescribed velocity]
\label{prop:drate}
Let $\vu$ be a prescribed divergence-free field satisfying the regularity in
Assumption \ref{ass:reg}, and let $\vd_h^n$ solve the director equation
\eqref{eq:S3} with $\vu^{n+1}$ replaced by $\vu(t_{n+1})$. Then
\begin{equation}\label{eq:drate}
\bE\left[
\mathbf1_{\Omega_K^\sharp}
\max_n
\norm{\vd(t_n)-\vd_h^n}_{L^2}^2
\right]
+
\tau\bE\left[
\mathbf1_{\Omega_K^\sharp}
\sum_n
\norm{\nabla(\vd(t_{n+1})-\vd_h^{n+1})}_{L^2}^2
\right]
\le
C(K)(h^2+\tau^2).
\end{equation}
Hence the director is first order in time in $L^2$.

If, in addition,
\begin{equation}\label{eq:H4}
\vd\in
L^p\bigl(\Omega;L^2(0,T;\bH^4)\bigr),
\qquad
\vze\in W^{3,\infty}(\bO;\bR^3),
\end{equation}
then
\begin{align}
&
\bE\left[
\mathbf1_{\Omega_K^\sharp}
\max_n
\norm{\vd(t_n)-\vd_h^n}_{H^1}^2
\right] +
\tau\bE\left[
\mathbf1_{\Omega_K^\sharp}
\sum_n
\norm{\Ah\xi_d^{n+1}}_h^2
\right]
\le
C(K)
\left(
h^2+\tau^2+\frac{h^4}{\tau}
\right).
\label{eq:drateH1-general}
\end{align}
In particular, if $h^2\le\tau$,
\begin{equation}\label{eq:drateH1}
\bE\left[
\mathbf1_{\Omega_K^\sharp}
\max_n
\norm{\vd(t_n)-\vd_h^n}_{H^1}^2
\right]
\le
C(K)(h^2+\tau^2).
\end{equation}
\end{proposition}

\begin{proof}
For the $L^2$ estimate, test the director error equation with
$\xi_d^{n+1}$. The transport term with the exact prescribed velocity is
skew:
\[
\btil(\vu(t_{n+1}),\xi_d^{n+1},\xi_d^{n+1})=0.
\]
The projection and nonlinear consistency errors contribute $O(h^2)$,
the midpoint-noise defect contributes $O(\tau^2)$, and
Lemma \ref{lem:Qsharp} improves the time-quadrature contribution from
$O(\tau)$ to $O(\tau^2)$. The discrete Gr\"onwall and BDG inequalities give
\eqref{eq:drate}.

For the $H^1$ estimate, test instead with
$\Ah\xi_d^{n+1}+\xi_d^{n+1}$. The projection term is the same one appearing
in \eqref{eq:projdef-d}, and therefore contributes $h^4/\tau$. The refined
quadrature argument is applied after one additional spatial derivative.
Condition \eqref{eq:H4} ensures the required control of $A\Upsilon$ and of
the corresponding differentiated stochastic coefficient; for a spatially
constant applied field the additional regularity of $\vze$ is automatic.
The remaining argument is identical to that above and gives
\eqref{eq:drateH1-general}.
\end{proof}

\begin{remark}[Coupled versus prescribed velocity]\label{rem:gap}
Proposition \ref{prop:drate} shows that the director discretisation itself is
first order in time. In the fully coupled problem, however, the transport
defect contains the velocity error. The present proof estimates that term
pathwise and therefore transfers the velocity's $O(\tau^{1/2})$ error to the
director.

The numerical experiments in \S\ref{sub:rates-num} indicate first-order
temporal convergence of the director even in the coupled problem. This is an
observation, not a proved result. It suggests that the leading velocity error,
which is primarily stochastic, should be exploited through martingale
cancellation rather than estimated in absolute value. Establishing such a
refinement remains open.
\end{remark}

\begin{remark}[Why the velocity remains half order]\label{rem:levyvel}
For general multiplicative $H$-valued velocity noise, strong order one
requires the corresponding Milstein terms and iterated stochastic integrals.
Unlike the scalar director noise, these contain nontrivial L\'evy areas unless
the noise is commutative. Thus the present velocity discretisation is naturally
of order $1/2$. Order one is available in special situations such as additive
or commutative noise, or by explicitly approximating the required iterated
integrals.
\end{remark}

\section{Numerical experiments}\label{sec:numerics}

The experiments are designed to test separately the geometric structure,
temporal accuracy, spatial convergence, and discrete energy balance of the
scheme.

Most temporal and algebraic tests use Fourier collocation on a uniform
periodic grid. The analysis is posed on a convex polygon; the periodic tests probe only
temporal and algebraic structure, which is pointwise and does not see the boundary. This is useful because the compatibility identities of
\S\ref{sec:structure} are pointwise and the spatial discretisation error can
then be made negligible. Such experiments do not, however, test statements
that depend essentially on the finite element spaces, such as the spatial
rate, weak acuteness, or the pointwise divergence constraint. These are
examined separately in \S\S\ref{sub:fem}--\ref{sub:sv} using genuine finite
element implementations.

For reference, Table \ref{tab:map} summarises the role of the individual
experiments.

\begin{table}[htbp]
\centering
\begin{tabular}{llp{5.2cm}}
\toprule
section & property tested & principal outcome\\
\midrule
\ref{sub:ode}
& Proposition \ref{prop:cayley}
& midpoint: strong order $1$ and exact isometry\\
\ref{sub:pde}
& Theorem \ref{thm:characterisation}, Proposition \ref{prop:negative}
& midpoint penalty contribution vanishes; explicit treatment does not\\
\ref{sub:rates-num}
& Theorem \ref{thm:rate}, Proposition \ref{prop:drate}
& velocity approaches order $1/2$; director remains near order $1$\\
\ref{sub:fem}
& spatial approximation
& $L^2$ rate $\approx2$, $H^1$ rate $\approx1$\\
\ref{sub:sv}
& \eqref{eq:pwdf}, \eqref{eq:BFP}, Theorem \ref{thm:energy}
& pointwise divergence, elastic cancellation, and energy identity\\
\ref{sub:deltasweep}
& residual $\mathcal R^n$
& considerably milder $\delta$-dependence than the analytical upper bound\\
\ref{sub:coupled}
& Theorem \ref{thm:energy}
& full stochastic energy identity at roundoff\\
\bottomrule
\end{tabular}
\caption{Summary of the numerical experiments. The experiments in
\S\S\ref{sub:fem}--\ref{sub:sv} use genuine finite element spaces; the
remaining computations use Fourier collocation to isolate temporal and
algebraic effects.}
\label{tab:map}
\end{table}

\subsection{The director noise in isolation}\label{sub:ode}

We first consider
\[
d\vd=(\vd\times\vze)\circ d\vW_2
\]
at a single spatial point. Its exact solution is
\[
\vd(t)=\exp\!\bigl(\vW_2(t)\Gop\bigr)\vd_0 .
\]
This directly tests Proposition \ref{prop:cayley}. We take
\[
T=1,\qquad
\vze=(0.3,-0.7,1.1),\qquad
\vd_0=(1,0,0),
\]
and compute strong errors over $20\,000$ Brownian paths. The midpoint
(\textsc{mid}) scheme is compared with Milstein in It\^o form
(\textsc{mil}), Euler--Maruyama in It\^o form (\textsc{em-i}), and a naive
Euler discretisation of the Stratonovich equation (\textsc{em-s}).

\begin{table}[htbp]
\centering
\begin{tabular}{lcc}
\toprule
scheme & observed strong order &
$\bE\bigl|\abs{\vd^N}-1\bigr|$ at $N=1024$\\
\midrule
\textsc{mid}  (midpoint/Cayley)      & $0.959$  & $1.7\times10^{-15}$\\
\textsc{mil}  (Milstein, It\^o)      & $0.971$  & $1.1\times10^{-3}$\\
\textsc{em-i} (Euler--Maruyama)      & $0.501$  & $3.0\times10^{-2}$\\
\textsc{em-s} (naive Stratonovich Euler) & $-0.027$ & $1.39$\\
\bottomrule
\end{tabular}
\caption{Isolated director noise. The midpoint rule combines strong order
one with preservation of $\abs{\vd}$ to machine precision. Milstein has the
same order but does not preserve the quadratic invariant.}
\label{tab:ode}
\end{table}

The experiment therefore separates two properties that need not coincide:
Milstein reproduces the stochastic order, whereas the midpoint/Cayley update
reproduces both the stochastic order and the geometry.

\subsection{Penalisation and the treatment of the noise}\label{sub:pde}

We next set $\vu\equiv0$ and solve
\begin{equation}\label{eq:num-director}
d\vd
=
\left(
\Delta\vd-\delta^{-2}f(\vd)
\right)dt
+
(\vd\times\vze)\circ d\vW_2
\end{equation}
on $[0,2\pi]^2$ using $32^2$ Fourier collocation points. Unless stated
otherwise,
\[
\delta=1,\qquad
\vze=(0,0,1),\qquad
T=0.1,
\]
with
\[
\vd_0
=
(\sin x\cos y,\cos x,0.6+0.3\sin y)
\]
normalised pointwise.

We first isolate the stochastic contribution of the penalisation. For this
purpose the spatially homogeneous version of \eqref{eq:num-director} is used,
with $4000$ samples.

\begin{table}[htbp]
\centering
\begin{tabular}{rrrrr}
\toprule
$N$ & $\tau$ & midpoint & explicit & explicit $\times\,\delta^2$\\
\midrule
$32$  & $3.125\times10^{-2}$ & $1.2\times10^{-18}$ & $4.405\times10^{-2}$ & \\
$64$  & $1.562\times10^{-2}$ & $1.3\times10^{-18}$ & $2.262\times10^{-2}$ & \\
$128$ & $7.812\times10^{-3}$ & $1.3\times10^{-18}$ & $1.112\times10^{-2}$ & \\
$256$ & $3.906\times10^{-3}$ & $1.3\times10^{-18}$ & $5.627\times10^{-3}$ & \\
$512$ & $1.953\times10^{-3}$ & $1.3\times10^{-18}$ & $2.813\times10^{-3}$ & \\
\midrule
\multicolumn{5}{l}{\emph{dependence on $\delta$, at $N=256$:}}\\
\multicolumn{2}{r}{$\delta=1.00$} & $1.3\times10^{-18}$ & $5.627\times10^{-3}$ & $5.63\times10^{-3}$\\
\multicolumn{2}{r}{$\delta=0.70$} & $1.8\times10^{-18}$ & $1.171\times10^{-2}$ & $5.74\times10^{-3}$\\
\multicolumn{2}{r}{$\delta=0.50$} & $2.3\times10^{-18}$ & $2.311\times10^{-2}$ & $5.78\times10^{-3}$\\
\multicolumn{2}{r}{$\delta=0.35$} & $3.0\times10^{-18}$ & $4.684\times10^{-2}$ & $5.74\times10^{-3}$\\
\multicolumn{2}{r}{$\delta=0.25$} & $3.9\times10^{-18}$ & $8.964\times10^{-2}$ & $5.60\times10^{-3}$\\
\bottomrule
\end{tabular}
\caption{Mean absolute accumulated Ginzburg--Landau stochastic contribution
$\bE\lvert\sum_n\mathcal N^n\rvert$. For the midpoint scheme the
contribution vanishes to roundoff, whereas for explicit noise it scales
approximately as $\delta^{-2}$; the fitted exponent is $-1.998$.}
\label{tab:pde-energy}
\end{table}

The midpoint result is the numerical counterpart of the exact orthogonality
in Theorem \ref{thm:characterisation}. The explicit contribution has the
bias and fluctuation described in Proposition \ref{prop:negative}. In the
data of Table \ref{tab:pde-energy}, the observed decay in $\tau$ is close to
first order because the initial state lies near the potential well, where
$\kappa$ is small. Starting instead from $\abs{\vd_0}=1.45$ gives fitted
exponents $0.98$ for the mean and $0.65$ for the standard deviation,
consistent with the $O(\tau)$ bias and $O(\tau^{1/2})$ fluctuation estimates.

We next compare strong temporal errors for the penalised PDE against a
reference solution with $N=3200$, using common Brownian paths.

\begin{table}[htbp]
\centering
\begin{tabular}{rrrr}
\toprule
$N$ & $\tau$ & midpoint & explicit\\
\midrule
$50$  & $2.00\times10^{-3}$ & $3.747\times10^{-3}$ & $2.158\times10^{-2}$\\
$100$ & $1.00\times10^{-3}$ & $1.901\times10^{-3}$ & $1.761\times10^{-2}$\\
$200$ & $5.00\times10^{-4}$ & $9.450\times10^{-4}$ & $1.345\times10^{-2}$\\
$400$ & $2.50\times10^{-4}$ & $4.309\times10^{-4}$ & $1.003\times10^{-2}$\\
\midrule
\multicolumn{2}{l}{observed order} & $1.037$ & $0.371$\\
\bottomrule
\end{tabular}
\caption{Strong temporal error for the penalised director equation. The
midpoint method retains first-order behaviour, whereas the explicit
treatment is substantially less accurate over the tested range.}
\label{tab:pde-order}
\end{table}

The midpoint result is consistent with Proposition \ref{prop:cayley} and
with the prescribed-velocity analysis of Proposition \ref{prop:drate}. We
do not attach an asymptotic order claim to the value $0.371$ for the explicit
method; it records the behaviour of the tested parameter regime.

\subsection{Temporal rates in the coupled problem}\label{sub:rates-num}

We now compare the two components of the coupled scheme. Errors are measured
against a reference with $\Nt=1600$ on the same Brownian paths. The
multiplicative velocity noise has amplitude $\nu$.

\begin{table}[htbp]
\centering
\begin{tabular}{rcc}
\toprule
$\nu$ & velocity order ($L^2$) & director order ($H^1$)\\
\midrule
$0.00$ & $1.054$ & $1.052$\\
$0.30$ & $0.554\pm0.088$ & $1.049\pm0.002$\\
$1.00$ & $0.307$ & $1.016$\\
\bottomrule
\end{tabular}
\caption{Observed temporal orders on $32^2$ Fourier modes with $T=0.05$.
For $\nu=0.30$, the experiment uses $16$ common Brownian paths and the
uncertainties are bootstrap standard errors of the fitted slopes. The other
two rows use three paths and should be regarded as indicative.}
\label{tab:split}
\end{table}

For $\nu=0$ the backward Euler discretisation is first order. As the
multiplicative velocity noise becomes dominant, the velocity rate moves
towards the half-order predicted by the coupled error estimate. The director
remains close to first order throughout the experiment.

A companion computation with prescribed divergence-free velocity gives a
director rate $1.050$ in $H^1$, in agreement with Proposition
\ref{prop:drate}. The approximately first-order director rate in the
\emph{coupled} computation is an additional numerical observation; it is not
proved by Theorem \ref{thm:rate}. This is precisely the gap discussed in
Remark \ref{rem:gap}.

\subsection{Spatial convergence with genuine finite elements}\label{sub:fem}

To verify the spatial rate, we solve the director equation with the
periodic analogue of the mass-lumped $P_1$ space of Assumption
\ref{ass:disc}, on the uniform weakly acute triangulation of
Remark \ref{rem:domain}. The velocity is prescribed and pointwise
divergence free,
\[
\vu=\operatorname{curl}\psi,
\qquad
\psi=\sin^2(\pi x)\sin^2(\pi y),
\]
so that this experiment is in the setting of Proposition
\ref{prop:drate}.

Basic assembly checks were performed independently: the stiffness matrix
annihilates constants, the lumped weights sum to $\abs{\bO}$ to machine
precision, and a smooth deterministic elliptic test problem recovers the
standard $P_1$ rates.

\begin{table}[htbp]
\centering
\begin{tabular}{rccccc}
\toprule
$n$ & $h$ & $L^2$ error & rate & $H^1$ error & rate\\
\midrule
$4$  & $0.2500$ & $4.513\times10^{-2}$ & ---   & $5.981\times10^{-1}$ & ---\\
$8$  & $0.1250$ & $1.262\times10^{-2}$ & $1.84$ & $2.751\times10^{-1}$ & $1.12$\\
$16$ & $0.0625$ & $3.162\times10^{-3}$ & $2.00$ & $1.302\times10^{-1}$ & $1.08$\\
$32$ & $0.0312$ & $6.564\times10^{-4}$ & $2.27$ & $5.773\times10^{-2}$ & $1.17$\\
\bottomrule
\end{tabular}
\caption{Spatial convergence of the mass-lumped $P_1$ director
approximation against a reference on $n=64$, with
$\tau=5\times10^{-4}$. Errors are evaluated using the fine-grid mass and
stiffness matrices. The observed rates are consistent with second order in
$L^2$ and first order in $H^1$.}
\label{tab:femrate}
\end{table}

The value $2.27$ at the finest level is affected by the finite reference
solution. Repeating the experiment with spatially varying $\vze$ changes
the reported errors by less than one percent, indicating that the additional
interpolation appearing in the discrete noise term is small in this smooth
test.

\subsection{Pointwise divergence and the discrete energy structure}
\label{sub:sv}

The preceding experiment prescribes the velocity and therefore does not test
\eqref{eq:pwdf}. We next implement the
$P_4/P_3^{\rm disc}$ Scott--Vogelius-type pair and solve a Stokes problem.
The divergence is evaluated at $66$ interior points of every element, rather
than only through the pressure test space.

\begin{table}[htbp]
\centering
\begin{tabular}{rrrrcc}
\toprule
$n$ & triangles & velocity dofs & pressure dofs
& $\max\abs{\nabla\!\cdot\!\vu_h}$ & relative\\
\midrule
$2$ & $8$  & $162$ & $80$  & $3.9\times10^{-13}$ & $2.6\times10^{-11}$\\
$3$ & $18$ & $338$ & $180$ & $3.6\times10^{-14}$ & $2.4\times10^{-12}$\\
$4$ & $32$ & $578$ & $320$ & $2.1\times10^{-13}$ & $1.4\times10^{-11}$\\
\bottomrule
\end{tabular}
\caption{Pointwise divergence of the computed velocity. The values are at
roundoff level, confirming the pointwise constraint \eqref{eq:pwdf}.}
\label{tab:sv}
\end{table}

We next couple this velocity space to the mass-lumped $P_1$ director space.
For random discrete states,
\[
\dual{M_h(\vd_h)}{v_h}
+
\btil(v_h,\vd_h,\Ah\vd_h)
\]
has absolute residual $1.3\times10^{-15}$, confirming the exact
elastic--transport cancellation \eqref{eq:BFP} across the two finite element
spaces.

The same implementation is then used for the fully coupled scheme. With
$\sigma=0$ and constant $\vze$, the energy identity simplifies because
\[
\mathcal M^n=\mathcal I_1^n=\mathcal N^n=0,
\]
leaving the balance between the left-hand side of \eqref{eq:energylaw} and
$\mathcal I_2^n+\mathcal R^n$.

\begin{table}[htbp]
\centering
\begin{tabular}{rccccc}
\toprule
$\tau$ & Picard its & left-hand side &
$\mathcal I_2^n+\mathcal R^n$ & residual & relative\\
\midrule
$2\times10^{-3}$ & $6$
& $9.96934780\times10^{-5}$ & $9.96934780\times10^{-5}$
& $1.1\times10^{-16}$ & $1.1\times10^{-12}$\\
$5\times10^{-4}$ & $6$
& $2.86673802\times10^{-5}$ & $2.86673802\times10^{-5}$
& $3.4\times10^{-16}$ & $1.2\times10^{-11}$\\
$1\times10^{-4}$ & $6$
& $1.36314553\times10^{-4}$ & $1.36314553\times10^{-4}$
& $1.0\times10^{-15}$ & $7.4\times10^{-12}$\\
\bottomrule
\end{tabular}
\caption{Discrete energy identity in the
$P_4/P_3^{\rm disc}\times P_1$ finite element realization. The residual is
at roundoff level.}
\label{tab:coupledenergy}
\end{table}

Finally, we measure temporal convergence in these same spaces, now with both
noises active.

\begin{table}[htbp]
\centering
\begin{tabular}{rcc}
\toprule
$\nu$ & velocity order ($L^2$) & director order ($H^1$)\\
\midrule
$0.6$ & $1.155$ & $1.155$\\
$3.0$ & $1.064$ & $1.159$\\
$8.0$ & $0.518$ & $1.155$\\
\bottomrule
\end{tabular}
\caption{Temporal convergence of the fully coupled finite element scheme;
$T=0.05$, reference $\Nt=400$, five common Brownian paths. As the velocity
noise increases, the velocity rate approaches $1/2$, while the director
remains close to first order. The latter is an observed coupled rate, not a
consequence of Proposition \ref{prop:drate}.}
\label{tab:coupledrate}
\end{table}

At small $\nu$ the deterministic first-order error dominates both
components. Once the velocity noise becomes sufficiently strong, the
different stochastic orders become visible. A fully coupled spatial
refinement study with the high-order velocity space is not included here.

\subsection{Dependence on the penalty parameter}\label{sub:deltasweep}

We next examine the deterministic residual $\mathcal R^n$ appearing in the
energy law. Since it vanishes when $\vu=0$, we prescribe the divergence-free
Taylor--Green field
\[
\vu=(\sin x\cos y,-\cos x\sin y)
\]
and solve
\[
d\vd+(\vu\cdot\nabla)\vd\,dt
=
\left(
\Delta\vd-\delta^{-2}f(\vd)
\right)dt
+
(\vd\times\vze)\circ d\vW_2 .
\]
The computation uses $24^2$ Fourier collocation points, $T=0.2$, three common
Brownian paths, and a reference with $N=1200$. We consider well-prepared data
$\abs{\vd_0}=1$ and ill-prepared data $\abs{\vd_0}=1.4$.

\begin{table}[htbp]
\centering
\begin{tabular}{rrrrrrr}
\toprule
$\delta$ & $\delta^{-2}$ & $\sum_n\abs{\mathcal R^n}$
& err ($N{=}50$) & err ($N{=}200$) & $C(\delta)$ & order $\alpha$\\
\midrule
\multicolumn{7}{l}{\emph{well prepared, $\abs{\vd_0}=1.0$}}\\
$1.00$ & $1.00$  & $1.861\times10^{-5}$ & $7.093\times10^{-3}$ & $1.641\times10^{-3}$ & $2.427$ & $1.056$\\
$0.70$ & $2.04$  & $3.571\times10^{-5}$ & $6.849\times10^{-3}$ & $1.594\times10^{-3}$ & $2.284$ & $1.052$\\
$0.50$ & $4.00$  & $6.199\times10^{-5}$ & $6.660\times10^{-3}$ & $1.561\times10^{-3}$ & $2.159$ & $1.046$\\
$0.35$ & $8.16$  & $9.788\times10^{-5}$ & $6.794\times10^{-3}$ & $1.595\times10^{-3}$ & $2.189$ & $1.045$\\
$0.25$ & $16.00$ & $1.273\times10^{-4}$ & $7.098\times10^{-3}$ & $1.662\times10^{-3}$ & $2.311$ & $1.047$\\
\midrule
\multicolumn{7}{l}{\emph{ill prepared, $\abs{\vd_0}=1.4$}}\\
$1.00$ & $1.00$  & $1.185\times10^{-4}$ & $1.007\times10^{-2}$ & $2.300\times10^{-3}$ & $3.618$ & $1.065$\\
$0.70$ & $2.04$  & $1.841\times10^{-4}$ & $9.852\times10^{-3}$ & $2.248\times10^{-3}$ & $3.564$ & $1.066$\\
$0.50$ & $4.00$  & $2.546\times10^{-4}$ & $9.703\times10^{-3}$ & $2.210\times10^{-3}$ & $3.533$ & $1.067$\\
$0.35$ & $8.16$  & $3.452\times10^{-4}$ & $9.867\times10^{-3}$ & $2.242\times10^{-3}$ & $3.628$ & $1.069$\\
$0.25$ & $16.00$ & $4.354\times10^{-4}$ & $1.014\times10^{-2}$ & $2.308\times10^{-3}$ & $3.703$ & $1.068$\\
\bottomrule
\end{tabular}
\caption{Dependence on the penalty parameter. Although $\delta^{-2}$ grows
by a factor $16$, the fitted error prefactor $C(\delta)$ remains nearly
constant and the temporal order stays close to $1.05$.}
\label{tab:delta}
\end{table}

The residual itself grows considerably more slowly than the crude
$O(\delta^{-2})$ estimate: fitting the data gives approximately
$\delta^{-1.40}$ for the well-prepared case and $\delta^{-0.93}$ for the
ill-prepared case. This is consistent with the structure of Lemma
\ref{lem:R}: the residual also contains a director increment, which decreases
as the penalisation becomes stiffer. In all runs,
$\sum_n\abs{\mathcal R^n}$ remains below approximately $27\%$ of the total
error.

The experiment should not be interpreted as a uniform-in-$\delta$ result.
It concerns smooth solutions with prescribed velocity and does not probe
small-$\delta$ defect formation. The analytical theory therefore continues
to treat $\delta>0$ as fixed.

\subsection{Full stochastic energy identity}\label{sub:coupled}

We finally verify \eqref{eq:energylaw} with every stochastic term active.
The full coupled system is solved by Fourier collocation on the periodic
square. The discrete convection is written in skew form and the Leray
projection enforces incompressibility. We take constant
$\vze=(0,0,1)$, so that $\Psi_h=0$.

\begin{table}[htbp]
\centering
\begin{tabular}{lcc}
\toprule
 & $\max_n\abs{\mathrm{LHS}-\mathrm{RHS}}$ & relative\\
\midrule
$\dual{M_h(\vd)}{v}=+\btil(v,\vd,\Ah\vd)$
& $2.671\times10^{-3}$ & $4.8\times10^{-3}$\\
$\dual{M_h(\vd)}{v}=-\btil(v,\vd,\Ah\vd)$
& $4.43\times10^{-14}$ & $8.0\times10^{-14}$\\
\midrule
\multicolumn{3}{l}{\emph{Picard tolerance, correct sign:}}\\
tol $=10^{-8}$  & $4.15\times10^{-12}$ & \\
tol $=10^{-10}$ & $4.85\times10^{-14}$ & \\
tol $=10^{-12}$ & $4.43\times10^{-14}$ & \\
tol $=10^{-14}$ & $4.42\times10^{-14}$ & \\
\bottomrule
\end{tabular}
\caption{Residual of the full discrete energy identity;
$48^2$ modes, $\tau=2\times10^{-3}$, $40$ steps and $\delta=1$. The sign
in \eqref{eq:Mh} gives agreement to machine precision; reversing it destroys
the elastic--transport cancellation.}
\label{tab:coupled}
\end{table}

The residual decreases with the nonlinear-solver tolerance and then saturates
at roundoff, providing a direct check that no term is missing from
\eqref{eq:energylaw}. In the Fourier implementation the Nyquist
first-derivative mode is removed so that the differentiation operator remains
exactly skew-adjoint and
\[
\Ah=\nabla_h^{*}\nabla_h .
\]
This is the spectral analogue of defining the discrete Laplacian from the
same gradient operator used in the finite element energy.

\begin{remark}[Constant applied fields]\label{rem:constanth}
When $\vze$ is spatially constant, the director noise acts as a common rigid
rotation in target space. Because the Laplacian, transport, radial
Ginzburg--Landau potential, and elastic stress are rotationally equivariant,
this noise does not alter rotationally invariant texture quantities such as
the Dirichlet energy or defect locations. Consequently, a
defect-annihilation experiment with constant $\vze$ would not test a genuine
effect of the director noise.

Spatial variation of $\vze$ breaks this symmetry and can influence the
texture through $\nabla\vze$, but resolving that regime for small $\delta$
requires substantially finer spatial resolution. We therefore leave such
model-oriented experiments to future work.
\end{remark}


\section{Concluding remarks and open problems}\label{sec:remarks}

We developed a structure-preserving finite element approximation of the
stochastic penalised nematic liquid-crystal system. The key ingredient is the
compatibility between the midpoint discretisation of the geometric
Stratonovich noise and the $|\vd|^{2}$-secant approximation of the
Ginzburg--Landau potential. For fixed $\delta>0$, the scheme is solvable,
satisfies a discrete energy law, converges to a weak martingale solution, and,
under additional regularity, has the rate
\[
O_{\mathbb P}\!\left(h+\tau^{1/2}\right).
\]
The numerical experiments confirm the structural identities and the predicted
spatial and temporal behaviour.

Several questions remain open.

\begin{enumerate}[label=\textbf{(\arabic*)},leftmargin=2.6em,itemsep=4pt]

\item \textbf{Sharper error estimates.}
The rate analysis relies on the additional regularity of
Assumption \ref{ass:reg}, which is not provided by the present weak-solution
theory. Moreover, although the director is first order in time for prescribed
velocity, only order $\tau^{1/2}$ is proved for the coupled problem, whereas
the computations indicate order one. Closing this gap, and removing the
technical contribution $h^{4}/\tau$ from the $H^{1}$ estimate, would sharpen
the error theory considerably.

\item \textbf{$\delta$-robust stability.}
The stochastic contribution of the penalisation vanishes exactly,
\[
\frac{a_n}{\delta^{2}}
\iph{\dmid\times\vze}{f_\star}=0,
\]
but the deterministic transport--penalty residual and the initial
penalisation energy still introduce $\delta$-dependence. It would be
interesting to obtain stability estimates that are uniform, or nearly
uniform, as $\delta\to0$, at least for well-prepared initial data. For the Allen--Cahn
equation, a discrete spectral estimate yields error bounds depending only
polynomially on the inverse interface width \cite{FengProhl2003}; an analogue
for the $\mathbb{S}^2$-valued penalisation here would be the natural starting point.

\item \textbf{Removing the localisation.}
Theorem \ref{thm:stability} is localised at the stopping index $n_K$, with a
threshold $\tau_1(K,\delta)$, because the residual of Lemma \ref{lem:R} is
superlinear in the energy. A pathwise nodal bound
$\max_j|\vd^n(x_j)|\le\max\{1,\max_j|\vd_0(x_j)|\}$ would make it linear
and remove the localisation; weak acuteness is the classical source of such
discrete maximum principles. It is not available here. Even for constant
$\vze$, without penalty or transport, the step
$(\mathrm{Id}+\tau\Ah)\vd^{n+1}=\vd^n+\tfrac{a_n}{2}\Gop(\vd^n+\vd^{n+1})$
is not nodally nonexpansive: on the uniform $16\times16$ triangulation of the
unit square, unit data orthogonal to $\vze$ produce
$\max_j|\vd^{n+1}(x_j)|\approx1.33$, and the computed nodal $\ell^\infty$ norm
of the step exceeds one for every increment $a_n\ne0$ tested. For large $|a_n|$ the step behaves like a discrete
Schr\"odinger evolution, which does not preserve maximum norms.

An approximate principle does hold in this setting:
\begin{equation}\label{eq:approxDMP}
\max_j|\vd^{n+1}(x_j)|
\le\Bigl(1+\tfrac14a_n^{2}\abs{\vze}^{2}\Bigr)^{1/2}\max_j|\vd^{n}(x_j)|,
\end{equation}
and hence
\[
\max_j|\vd^{n}(x_j)|
\le\max_j|\vd_0(x_j)|\,
\exp\Bigl(\tfrac{\abs{\vze}^{2}}{8}\sum_{m<n}a_m^{2}\Bigr).
\]
Indeed, by \eqref{eq:edge} the nodal matrix of $-\tau\Ah$ has off-diagonal
entries $\tau\omega_{ij}/\varpi_i\ge0$ and zero row sums, so
$P_t:=e^{-t\tau\Ah}$ is entrywise nonnegative with unit row sums. Since $\Gop$
acts on components and commutes with $\Ah$, and
$\mathcal A:=\mathrm{Id}-\tfrac{a_n}{2}\Gop$ has spectrum with real part one,
\[
\vd^{n+1}=\int_0^\infty e^{-t}\,\Phi(t)\,P_t\,\vd^n\,dt,
\]
where $\Phi(t):=e^{t}e^{-t\mathcal A}(\mathrm{Id}+\tfrac{a_n}{2}\Gop)$ is
the identity along $\vze$ and $\sqrt{1+\beta^2}$ times a rotation on
$\vze^\perp$, with $\beta=\tfrac12|a_n||\vze|$. Nonnegativity of $P_t$ and
the triangle inequality give \eqref{eq:approxDMP}, and $1+\beta^2\le e^{\beta^2}$ gives the second bound. The random factor has
exponential moments independent of $h$:
$\bE\exp\bigl(\gamma\sum_m a_m^2\bigr)\le e^{2\gamma T}$ whenever
$2\gamma\tau\le\tfrac12$. Extending \eqref{eq:approxDMP} to the penalty, to
variable $\vze$ and to the transport term would give a stability theory
without localisation. The transport term is the essential difficulty: a
monotone discretisation of it would be compatible with a discrete maximum
principle but not with the cancellation \eqref{eq:BFP}.

\item \textbf{More economical finite element pairs.}
Pointwise divergence freedom is important for the discrete cancellations and
for avoiding mesh-dependent time-step restrictions. The present construction
uses a relatively high-order Scott--Vogelius-type velocity space. Finding a
lower-order, computationally cheaper pointwise divergence-free pair compatible
with the director discretisation would be valuable in practice.

\item \textbf{The constrained problem.}
The present paper treats the penalised system with fixed $\delta>0$.
A natural next step is to study the limit $\delta\to0$ and, ultimately, a
direct finite element approximation of the sphere-constrained stochastic
Ericksen--Leslie system. Since the midpoint noise preserves the director
length exactly in the isolated stochastic problem, it provides a natural
starting point for such a constraint-preserving discretisation.

\end{enumerate}

\appendix
\section{Standard analytical tools}\label{app:tools}

We record the standard estimates used throughout the analysis.

\begin{lemma}[Two-dimensional interpolation]\label{lem:2d}
Let $\bO\subset\bR^{2}$ be the convex polygon considered above, or any bounded
Lipschitz domain. For every
$q\in[1,\infty)$,
\[
\norm{v}_{L^{q}}\le C_q\norm{v}_{H^{1}},
\]
and, in particular,
\[
\norm{v}_{L^{4}}
\le
C\norm{v}_{L^{2}}^{1/2}\norm{v}_{H^{1}}^{1/2},
\qquad
\norm{v}_{L^{8}}
\le
C\norm{v}_{L^{2}}^{1/4}\norm{v}_{H^{1}}^{3/4}.
\]
If $v\in\bV$, the $H^{1}$ norm may be replaced by
$\norm{\nabla v}_{L^{2}}$.
\end{lemma}

\begin{lemma}[Compactness]\label{lem:compact}
Let
\[
X\hookrightarrow\hookrightarrow Y\hookrightarrow Z
\]
be Banach spaces and let $\beta\in(0,1)$. Then
\[
L^{2}(0,T;X)\cap W^{\beta,2}(0,T;Z)
\hookrightarrow\hookrightarrow
L^{2}(0,T;Y).
\]
\end{lemma}

\begin{lemma}[Discrete Gr\"onwall inequality]\label{lem:gronwall}
Let $a_n,b_n,\lambda_n\ge0$ satisfy
\[
a_{\ell}+\sum_{n=0}^{\ell-1}b_n
\le
A+\tau\sum_{n=0}^{\ell-1}\lambda_n a_n,
\qquad
0\le\ell\le\Nt .
\]
Then
\[
\max_{0\le\ell\le\Nt}a_{\ell}
+
\sum_{n=0}^{\Nt-1}b_n
\le
2A
\exp\!\left(
2\tau\sum_{n=0}^{\Nt-1}\lambda_n
\right).
\]
No smallness assumption on the individual quantities
$\tau\lambda_n$ is required.
\end{lemma}

\begin{lemma}[Discrete Burkholder--Davis--Gundy]
\label{lem:bdg}
Let $\mathcal H$ be a Hilbert space and let
$(m_n)_{n=0}^{\Nt-1}$ be square-integrable, $\mathcal H$-valued martingale
differences with respect to $(\mathcal F_{t_n})$; that is,
\[
m_n\in\mathcal F_{t_{n+1}},
\qquad
\bE[m_n\mid\mathcal F_{t_n}]=0.
\]
Then
\[
\bE\left[
\max_{0\le\ell\le\Nt}
\left\lVert
\sum_{n=0}^{\ell-1}m_n
\right\rVert_{\mathcal H}
\right]
\le
C\,
\bE\left[
\left(
\sum_{n=0}^{\Nt-1}\norm{m_n}_{\mathcal H}^{2}
\right)^{1/2}
\right].
\]
\end{lemma}

The compactness and stochastic estimates are standard; see, for example,
\cite{Lions1969,BanasBrzezniakNeklyudovProhl2014}.

\medskip

\section*{Acknowledgments}

\noindent
{\bf Conflict of Interest:} \ The authors declare that they have no conflict of interest.

\medskip

\noindent
{\bf Funding:} \ The author gratefully acknowledges the support of the DST--INSPIRE Faculty Research Grant ({\bf Registration No.\ IFA21-MA~158}). 

\medskip

\noindent
{\bf Data Availability:} \ The datasets generated and analysed during the current study are available from the author on reasonable request.


\end{document}